\documentclass{conm-p-l}
\usepackage{amsmath,amscd,amssymb,amsfonts}
\usepackage{hyperref}
\copyrightinfo{2002}{American Mathematical Society}

\newtheorem{thm}{Theorem}[section]
\newtheorem{lm}{Lemma}[section]

\newtheorem{prop}{Proposition}[section]
\newtheorem{cor}{Corollary}[section]

\theoremstyle{definition}
\newtheorem{de}{Definition}[section]
\newtheorem{ex}{Example}[section]
\newtheorem*{exs}{Examples}

\theoremstyle{remark}
\newtheorem{rem}{Remark}[section]

\numberwithin{equation}{section}

\def\co{\colon\thinspace}

\newcommand{\hash}{\#}

\renewcommand{\div}{\mathrm{div\,}}
 
\DeclareMathOperator{\tr}{tr} \DeclareMathOperator{\diag}{diag}
\DeclareMathOperator{\antidiag}{antidiag} \DeclareMathOperator{\str}{str}
\DeclareMathOperator{\otr}{otr}

\DeclareMathOperator{\Der}{Der} \DeclareMathOperator{\Vect}{Vect}
\DeclareMathOperator{\ad}{ad}
\newcommand{\der}[2]{{\frac{\partial {#1}}{\partial {#2}}}}
\newcommand{\lder}[2]{{\partial {#1}/\partial {#2}}}
\newcommand{\dder}[3]{{\frac{\partial^2 {#1}}{\partial {#2}\partial {#3}}}}
\newcommand{\R}[1]{{\mathbb R}^{#1}}
\newcommand{\RR}{\mathbb R}
\newcommand{\Z}{{\mathbb Z_{2}}}
\newcommand{\ZZ}{{\mathbb Z}}

\renewcommand{\d}{\partial}

\newcommand{\fun}{C^{\infty}}
\newcommand{\widebar}{\overline}

\newcommand{\G}{\Gamma}
\renewcommand{\a}{\alpha}

\newcommand{\dd}{\delta}
\newcommand{\e}{\varepsilon}

\newcommand{\h}{\eta}
\renewcommand{\k}{\varkappa}
\renewcommand{\t}{\theta}
\renewcommand{\l}{\lambda}

\newcommand{\x}{\xi}
\renewcommand{\o}{\omega}

\newcommand{\p}{\pi}

\newcommand{\kt}{{\tilde k}}
\newcommand{\jt}{{\tilde \jmath}}
\newcommand{\itt}{{\tilde \imath}}
\newcommand{\lt}{{\tilde l}}
\newcommand{\at}{{\tilde a}}
\newcommand{\bt}{{\tilde b}}
\newcommand{\ct}{{\tilde c}}
\newcommand{\dt}{{\tilde d}}
\newcommand{\ft}{{\tilde f}}
\newcommand{\gtt}{{\tilde g}}

\newcommand{\ut}{{\tilde u}}
\newcommand{\vt}{{\tilde v}}
\newcommand{\xt}{{\tilde x}}
\newcommand{\Xt}{{\tilde X}}
\newcommand{\Qt}{{\tilde Q}}
\newcommand{\Pt}{{\tilde P}}

\begin{document}

\title{Graded manifolds and Drinfeld doubles for Lie
bialgebroids}

\author{Theodore Voronov}
\address{Department of Mathematics, University of Manchester Institute of Science
and Technology (UMIST), PO Box 88, Manchester M60 1QD, England}
\email{theodore.voronov@umist.ac.uk}
\thanks{The research was supported in part by EPSRC, under grant GR/N00821/01.}


\date{}

\begin{abstract}
\noindent We define \textit{graded manifolds} as a version of
supermanifolds endowed with an extra $\mathbb Z$-grading in the
structure sheaf, called \textit{weight} (not linked with parity).
Examples are ordinary supermanifolds, vector bundles, double
vector bundles (in particular, iterated constructions like $TTM$),
etc. I give a construction of \textit{doubles} for \textit{graded}
$QS$- and \textit{graded $QP$-manifolds} (graded manifolds endowed
with a homological vector field and a Schouten/Poisson bracket).
Relation is explained with Drinfeld's Lie bialgebras and their
doubles. Graded $QS$-manifolds can be considered, roughly, as
``generalized  Lie bialgebroids''. The double for them is closely
related with the analog of Drinfeld's double for Lie bialgebroids
recently suggested by Roytenberg.  Lie bialgebroids as a
generalization of Lie bialgebras, over some base manifold, were
defined by Mackenzie and P.~Xu. Graded $QP$-manifolds give an {odd
version} for all this, in particular, they contain  ``odd
analogs'' for Lie bialgebras, Manin triples, and Drinfeld's
double.
\end{abstract}

\maketitle

\section*{Introduction}
Analog for Lie bialgebroids of Drinfeld's classical  double is a
puzzle. Various constructions were suggested (Mackenzie,
Liu--Weinstein--Xu and Roytenberg), all different and all giving
something which is not a Lie bialgebroid, even not a Lie
algebroid. Recall that the Drinfeld double of a Lie bialgebra is
again a Lie bialgebra (with nice properties). In Roytenberg's
thesis~\cite{roytenberg:thesis} it was suggested to consider the
supermanifold $T^*(\Pi E)$ (with a particular natural structure)
as the ``right'' notion of the double for a Lie bialgebroid $E\to
X$ over a base $X$. (Close ideas are due to
A.~Vaintrob~\cite{vaintrob:algebroids}.)

In this paper we define a class of supermanifolds which, in
particular, generalizes the structure of Lie bialgebroids. We
prove that an analog of Roytenberg's construction carries over for
them. Namely, we consider a ``graded manifold'' (see below) $M$
endowed with a homological vector field $\hat Q$ of weight $q$ and
a compatible Schouten bracket of weight $s$, and  prove that on
the space of the cotangent bundle $T^*M$ there is a naturally
defined homological field $\hat Q_D$ and it is possible to
introduce weight for $T^*M$ in such a way that $\hat Q_D$ will be
again of the weight $q$ (same as original $\hat Q$), and certain
natural properties hold. The space $DM=T^*M$ with such a structure
is called the \textit{double} of $M$. The double $DM$ so defined
inherits half the original structure of $M$, a homological field.
Using a linear connection on $M$, it is possible to define on $DM$
an ``almost'' Schouten bracket as well. I show that this exactly
recovers the whole structure of the Drinfeld double when
restricted to Lie bialgebras. (For a nonlinear $M$, the bracket on
$DM$ does not, in general, obey the Jacobi identity.) By
considering supermanifolds with a Poisson rather than Schouten
bracket, we obtain an {odd  analog} of this theory. We introduce
and study \textit{odd Lie bialgebras} and  \textit{odd double} for
them. In particular, we construct an odd Lie bialgebra structure
in the Lie superalgebra $\mathfrak{q}(n)$. (A quantization of this
example should give a quantum supergroup $GQ_{\e}(n)$ with an odd
quantum parameter $\e$.)

This work is based on two simple ideas.

The first is a key notion of a {graded manifold}. Probably, the
concept is ``in the air''. Needless to argue for the usefulness of
$\mathbb Z$-grading in algebraic questions. Unfortunately, there
is an often confusion with $\Z$-grading  responsible for the sign
rule. Of course, in many examples the $\Z$-grading is induced by
some underlying $\mathbb Z$-grading, but in other examples not,
and a natural $\mathbb Z$-grading can have nothing to do with the
sign rule (example:  degree for polynomials). We introduce
\textit{graded manifolds} as a certain version of supermanifolds
with a $\mathbb Z$-grading in the structure sheaf
\textbf{independent of parity}. We call it \textit{weight}. This
structure embraces many examples and carries over to natural
constructions. (The term ``graded manifolds'' was used in some
early works on supermanifolds instead of the latter term. Our
usage has no relation with this.)

The second is a ``geometrization'' of the concept of a Lie
algebra, in the following sense. (Everywhere when we say ``Lie
algebra'' we  mean ``Lie superalgebra''.) This is the answer to a
simple question: to which structure in the algebra of functions
corresponds the Lie bracket? Consider simultaneously the vector
space $\mathfrak g$ with all its ``neighbors'': $\mathfrak g^*$,
$\Pi \mathfrak g$, $\Pi \mathfrak g^*$. On $\mathfrak g^*$ we get
a linear Poisson bracket (the Lie--Berezin--Kirillov bracket). On
$\Pi \mathfrak g^*$ we get a linear odd Poisson (=Schouten,
Gerstenhaber) bracket, which we call the Lie--Schouten bracket. On
$\Pi \mathfrak g$ we get a quadratic homological vector field
$Q=\pm (1/2) c_{ij}^k\x^j\x^i \lder{}{\xi^k}$. Hence, Poisson
brackets, Schouten  brackets and homological vector fields  are
all, equally, extensions of a Lie algebra structure. We call
supermanifolds with such structures \textit{$P$-manifolds},
\textit{$S$-manifolds} and \textit{$Q$-manifolds},
respectively\,\footnote{Curiously, the most difficult is to tell
what structure in the algebra of functions on $\mathfrak g$ itself
corresponds to the Lie bracket on $\mathfrak g$. An analogous
question is, what natural operator exists on functions on $TM$,
i.e., on Lagrangians. Partial answer, of course, is: the
Euler--Lagrange equations.  Their analog for Lie algebras and Lie
algebroids was considered by
A.~Weinstein~\cite{weinstein:lagrangian}.}. The next question is,
what is a ``bi-'' structure from this viewpoint. Shortly, it is a
Lie algebra structure on two neighbors with a compatibility
condition. We get the following list: $(\mathfrak g,\mathfrak
g^*)$, which is a Lie bialgebra in Drinfeld's sense, $(\mathfrak
g,\Pi \mathfrak g^*)$, which is the {odd analog} of Drinfeld's
notion (can be non-trivial only for superalgebras), and
$(\mathfrak g,\Pi \mathfrak g)$. Their geometric versions are,
respectively: \textit{$QS$-manifolds} (compatible $Q$- and
$S$-structures), \textit{$QP$-manifolds} (compatible $Q$- and
$P$-structures), and \textit{$PS$-manifolds} (two compatible
brackets of the opposite parity;  structures of this kind have
very interesting geometry, their study was initiated
in~\cite{hov:delta}). Notice, that for a ``bi-'' Lie structure of
each type there is just one geometric manifestation (one type of
structure, though it can be realized on different spaces, e.g.,
the $QS$-structures on $\Pi \mathfrak g$ and  $\Pi \mathfrak g^*$
for a Drinfeld Lie bialgebra $\mathfrak g$), compared to the three
different manifestations of a ``single'' Lie bracket.

Combining these two ideas, we say that instead of conventional Lie
(bi)\-algeb\-roids one should deal with graded manifolds with
$QS$- and $QP$-structure and keep track of weights. In particular,
this allows to develop a nice framework for doubles. As a
byproduct, we naturally encounter odd analogs for Lie bialgebras
and bialgebroids.

The paper is organized as follows:

In Section~\ref{psalgebras} we review the necessary facts concerning Poisson and Schouten
algebras. In particular, it is shown how an arbitrary Poisson structure is expressed via
the canonical Schouten structure on the anticotangent bundle, and an arbitrary Schouten
structure is expressed via the canonical Poisson bracket on the cotangent bundle. We see
how even and odd brackets intertwine.

In Section~\ref{bistructures} we give a geometric description of
Drinfeld's Lie bialgebras, in the language of supermanifolds, and
we explain how the classical Drinfeld's double for $\mathfrak g$
is equivalent to a natural construction of a $QS$-structure on
$T^*\Pi\mathfrak g\cong \Pi(\mathfrak g\oplus\mathfrak g^*)$. This
also leads us to a new notion of an \textit{odd Lie bialgebra},
which we define in terms of a $QP$-structure on $\Pi \mathfrak g$.

In Section~\ref{bialgebroids} we review the ``super'' approach to Lie bialgebroids and
doubles. It is a generalization of the ``super'' description of Lie bialgebras and their
Drinfeld doubles given in Section~\ref{bistructures}. We describe Roytenberg's ``Drinfeld
double'' of a Lie bialgebroid. It is no longer a Lie algebroid. Actually, it is an
example of a graded manifold (which we define in the next section).

In Section~\ref{graded} we introduce graded manifolds as
supermanifolds with a ``weight'' in the algebra of functions, and
show how it carries over to natural bundles and the canonical
brackets. In particular, we show that the natural bundles like
$TM$ are bi-graded. This gives important flexibility for the
definition of a ``total'' weight.

In Section~\ref{generalized} we consider graded $QS$- and  $QP$-manifolds. The former
should be regarded as generalized Lie (anti)bialgebroids. The latter are a generalization
of odd Lie bialgebras defined in Section~\ref{bistructures}. For both we give a
construction of a \textit{double} as a graded $Q$-manifold. For a graded $QS$-manifold
$\Pi E$, where $E$ is a Lie bialgebroid, this recovers Roytenberg's construction. We show
how with the help of a connection it is possible to construct an ``almost'' bracket in
the double.

In Section~\ref{odd} we consider {odd Lie bialgebras} and the \textit{odd double} in
detail. We consider as example the matrix superalgebra $\mathfrak q(n)$ and  show that it
has a natural odd bialgebra structure. As a technical tool to define this structure, we
introduce a ``relative'' version of doubles, in the situation when there is an action of
some Lie superalgebra.  (This models certain triangular decompositions of superalgebras.)

In Appendix we give a proof of a theorem of Mackenzie and Xu about
cotangent bundles for dual vector bundles. We also suggest and
prove its odd analog.

\textbf{Note about usage.} All our constructions are ``superized'', i.e., we
automatically work with supermanifolds, superalgebras, super Lie (bi)algebroids, etc.,
from the beginning.  ``Algebras'' and ``manifolds'' stand for ``superalgebras'' and
``supermanifolds'' throughout the text unless otherwise required for clarity. $\Pi$
stands for the parity reversion functor. We use the prefix ``anti'' for \textit{opposite}
vector spaces and vector bundles (spaces and bundles with reversed parity).  For example,
\textit{anticotangent bundle} means $\Pi T^*M$,  \textit{Lie antialgebra} means the space
$\Pi \mathfrak g$ for a Lie algebra $\mathfrak g$,  an \textit{antivector} (for a vector
space $V$) means an element of $\Pi V$.

\section{Poisson and Schouten algebras}\label{psalgebras}

For the sake of completeness, we provide some definitions and
facts. All vector spaces are assumed $\Z$-graded, with tilde
denoting parity of a homogeneous element. Commutators,
commutativity, derivations, linearity, etc., are understood in the
$\Z$-graded sense. The numbers $0$, $1$ denoting the grading are
residues modulo $2$.

\begin{de}
A \textit{Poisson algebra} is a vector space $A=A_0\oplus A_1$ with a bilinear associative
multiplication $(a,b)\mapsto ab$ such that $A_iA_j\subset A_{i+j}$ and an even bilinear operation
denoted $\{\ ,\ \}$ (i.e., $\{A_i,A_j\}\subset A_{i+j}$) satisfying
\begin{align}
  \{a,b\}&=-(-1)^{\at\bt}\{b,a\}    \label{panti}\\
  \{a,\{b,c\}\} & =\{\{a,b\},c\}+(-1)^{\at\bt}\{b,\{a,c\}\} \\
  \{a,bc\} & =\{a,b\}c+(-1)^{\at\bt}b\{a,c\} \label{poisson}
\end{align}
for every homogeneous elements $a,b,c\in A$. The operation $\{\ ,\ \}$ is called \textit{Poisson
bracket}.
\end{de}

\begin{de}
A \textit{Schouten algebra} is a vector space $A=A_0\oplus A_1$ with a bilinear associative
multiplication $(a,b)\mapsto ab$ such that $A_iA_j\subset A_{i+j}$ and an odd bilinear operation
denoted $\{\ ,\ \}$ (i.e., $\{A_i,A_j\}\subset A_{i+j+1}$) satisfying
\begin{align}
  \{a,b\}&=-(-1)^{(\at+1)(\bt+1)}\{b,a\} \label{santi} \\
  \{a,\{b,c\}\} & =\{\{a,b\},c\}+(-1)^{(\at+1)(\bt+1)}\{b,\{a,c\}\} \\
  \{a,bc\} & =\{a,b\}c+(-1)^{(\at+1)\bt}b\{a,c\} \label{schouten}
\end{align}
for every homogeneous elements $a,b,c\in A$. The operation $\{\ ,\ \}$ is called \textit{Schouten
bracket}.
\end{de}

The difference between these two notions is obvious. A Poisson
algebra is, in particular, a Lie (super)algebra, and the
identity~\eqref{poisson} tells that for every $a\in A$, the
operator $\ad a=\{a,\ \}$ acts as a derivation of the associative
product. The operator $\ad a$ has the same parity as $a$.
Similarly, a Schouten algebra is, in particular,  a Lie
superalgebra, but  w.r.t. the reversed parity. The operator $\ad
a:=(-1)^{\at+1}\{a,\ \}$ in a Schouten algebra is a derivation of
the associative product of parity opposite to that of $a$. No
shift of parity can turn~\eqref{schouten} into~\eqref{poisson}.

Other names for Schouten algebras are \textit{odd Poisson
algebras} and \textit{Gerstenhaber algebras}. Schouten bracket is
also known as \textit{odd bracket} and \textit{antibracket}.

\begin{rem}
Schouten algebras became very popular recently, in connection with
their use in Batalin--Vilkovisky (BV) formalism in quantum field
theory~\cite{bv:perv},\cite{bv:vtor} and in deformation
quantization (see, e.g.,~\cite{kontsevich:quant}). The name
``Gerstenhaber algebras'' is motivated by the corresponding
structure discovered by
Gerstenhaber~\cite{gerstenhaber:cohomology63} in the cohomology of
an associative algebra. The name ``Schouten algebra'' comes from
the canonical Schouten bracket of multivector
fields~\cite{schouten:1940}, see below. This bracket is completely
parallel with the canonical Poisson bracket on $T^*M$, so  the
choice of ``Schouten algebras'' for the odd counterpart of
``Poisson algebras'' seems natural. For a manifold provided with a
volume element there is another operation on multivector fields,
the divergence. It satisfies the relation
\begin{equation}\label{divergence}
  \dd(PQ)=\dd P \,Q +(-1)^{\Pt}P\,\dd Q +(-1)^{\Pt +1}\{P,Q\},
\end{equation}
thus connecting Schouten bracket and multiplication. Schouten algebras endowed with a
differential satisfying~\eqref{divergence} received the name \textit{Batalin--Vilkovisky
algebras} (BV--algebras), because of their role in  the BV--formalism.
\end{rem}

\begin{de}A (super)manifold with a Poisson (Schouten) structure in the algebra of functions is
called a \textit{Poisson} (resp., \textit{Schouten} or \textit{odd Poisson}) \textit{manifold}.
\end{de}

Basic examples: $T^*M$ and the canonical Poisson bracket on
functions on $T^*M$ (Hamiltonians on $M$), $\Pi T^*M$ and the
canonical Schouten bracket on functions on $\Pi T^*M$ (multivector
fields on $M$). They solve two universal problems: for Lie
homomorphisms $\Vect M\to A$ to commutative Poisson algebras
``over $C^{\infty}(M)$'' and for odd Lie homomorphisms $\Vect M\to
A$ to commutative Schouten algebras ``over $C^{\infty}(M)$''
(precise formulation requires the notions of Lie--Rinehart
algebras or Lie algebroids, see below). In the sequel we will use
the following notation for the homomorphisms from vector fields to
Hamiltonians and to multivector fields: $p{\co} X\mapsto p_X$ and
$\theta{\co}  X\mapsto \theta_X$, respectively. (Notice that the
second map is odd.)

A purely algebraic variant is given by the Poisson bracket on
$\fun(\mathfrak g^*)$,
\begin{equation}\label{liep}
  \{x_i,x_j\}=c_{ij}^kx_k,
\end{equation}
for a Lie algebra $\mathfrak g$, the
\textit{Lie--Poisson(--Berezin--Kirillov) bracket}, and by the
Schouten bracket on $\fun(\Pi \mathfrak g^*)$,
\begin{equation}\label{lies}
  \{\x_i,\x_j\}=c_{ij}^k\x_k,
\end{equation}
(where the coordinates $\x_i$ have parity opposite to that of the
respective $x_i$), the \textit{Lie--Schouten bracket}, as we shall
call it. They solve universal problems for even/odd
bracket-preserving maps $\mathfrak g\to A$ to commutative
Poisson/Schouten algebras respectively. (For the even case this is
directly related with the moment map for Hamiltonian actions of
$\mathfrak g$.)

Examples above are united in the concept of a Lie algebroid
(see~\cite{mackenzie:book,mackenzie:alg}). A \textit{Lie
algebroid} is a vector bundle $E\to M$ over a manifold $M$ endowed
with a vector bundle map $a{\co} E\to TM$ (called the
\textit{anchor}) and a Lie algebra structure on the space of
sections $\fun(M,E)$ satisfying
\begin{align}
  [u,fv]&={a(u)}f\,v+(-1)^{\ut\ft}f[u,v] \label{lr1}\\
  a([u,v])&= [a(u),a(v)]\label{lr2}
\end{align}
for all $u,v\in \fun(M,E)$ and all $f\in \fun(M)$. Here we
identify vector fields on $M$ with derivations of the algebra
$\fun(M)$. Likewise, for an arbitrary (commutative associative)
algebra $A$,  a \textit{Lie  pseudoalgebra} or a
\textit{Lie--Rinehart algebra} over $A$ is a module $L$ over $A$
with a Lie algebra structure and a homomorphism of modules $a{\co}
L\to \Der A$ satisfying~\eqref{lr1},\eqref{lr2}
(see~\cite{mackenzie:alg}).

It is not difficult to see that every Lie algebroid (or every
Lie--Rinehart algebra) canonically produces two algebras, a
Poisson algebra and a Schouten algebra. The Poisson algebra
associated with a Lie algebroid  $E$ is $\fun(E^*)$. The Schouten
algebra associated with   $E$ is $\fun(\Pi E^*)$. Axioms of
Poisson/Schouten algebra for $\fun(E^*)$, $\fun(\Pi E^*)$ are
equivalent to $E$ being a Lie algebroid. These algebras are
universal in a certain precise sense, and their construction
unites the constructions of the canonical Poisson bracket in
$\fun(T^*M)$ and the canonical Schouten bracket in $\fun(\Pi
T^*M)$ from the commutator of vector fields with the constructions
of the Lie--Poisson bracket in $\fun(\mathfrak g^*)$ and the
Lie--Schouten bracket in $\fun(\Pi \mathfrak g^*)$ from the
bracket in $\mathfrak g$.

\smallskip
{\footnotesize Most often the associative multiplication in the
definition of Poisson/Schouten algebras is assumed to be
commutative.  There are natural examples where it is not so. The
commutator bracket $[a,b]=ab-(-1)^{\at\bt}ba$ obviously satisfies
the Leibniz identity~\eqref{poisson} and thus makes every
associative algebra into a (noncommutative) Poisson algebra. Let
us reserve the notation $[\ ,\ ]$ for the commutator bracket and
denote a given arbitrary Poisson bracket by $\{\ ,\ \}$.

\begin{lm}[\cite{tv:poiss}]
In every Poisson algebra the following identity holds:
\begin{equation}\label{myidentity}
  [a,b]\{c,d\}=\{a,b\}[c,d],
\end{equation}
for arbitrary elements $a,b,c,d$.
\end{lm}
Obviously, \eqref{myidentity} degenerates in two special cases,
when the algebra is commutative and when the bracket coincides
with the commutator, as in the example above.

Consider a Lie algebra $\mathfrak g$. Let us denote the Lie brackets in $\mathfrak g$ by $\{\ ,\
\}$. One may ask about bracket-preserving maps $\mathfrak g\to A$ to Poisson algebras without
conditions of commutativity.
\begin{thm}[\cite{tv:poiss}] A non-commutative Poisson algebra that is a universal range for such
maps (the ``Poisson envelope'' of  a Lie algebra $\mathfrak g$) is given by $E(\mathfrak
g):=T(\mathfrak g)/I$ where the ideal $I$ in the free tensor algebra is generated by the
elements of the form
\begin{equation}\label{ideal}
  \bigl(a\otimes b-(-1)^{\at\bt}b\otimes a\bigr)\otimes\{c,d\}-
  \{a,b\}\otimes\bigl(c\otimes d-(-1)^{\ct\dt}d\otimes c\bigr),
\end{equation}
for all $a,b,c,d\in \mathfrak g$.
\end{thm}
The algebra $\fun_{pol}(\mathfrak g^*)$ and the universal
enveloping algebra $U(\mathfrak g)$ both are quotients of
$E(\mathfrak g)$. The algebra $E(\mathfrak g)$ can be viewed as a
non-commutative analog of the Kirillov--Kostant--Sourieau moment
space $\mathfrak g^*$. For details see~\cite{tv:poiss}.

}

\smallskip

On a manifold $M$, a non-degenerate bracket is specified by a
symplectic $2$-form $\omega\in\Omega^2(M)$,  even for the Poisson
case and odd for the Schouten case. The correspondence is given by
the formulas: $i_{X_f}\omega=-(-1)^{s\ft}df$,
$\{f,g\}=(-1)^{s\ft}X_fg$, $X_f$ being the Hamiltonian vector
field corresponding to a function $f$, where $s=0$ for an even
bracket and $s=1$  for an odd bracket.

Consider the general situation.
From~\eqref{poisson},\eqref{panti},\eqref{schouten},\eqref{santi}
follows that a Poisson as well as a Schouten bracket is a
bi-derivation of functions. Hence a bracket is specified by a
contravariant tensor field of rank $2$. From antisymmetry follows
that in the Poisson case it should be an even bivector field. What
about the Schouten case?  Contrary to what one might think at the
first moment, an odd bracket is not related with  odd bivector
fields. In fact, a Schouten bracket on $M$ is specified by a
fiberwise quadratic odd Hamiltonian function\footnote{In a
non-degenerate situation, brackets of both types are specified by
$2$-forms. Noteworthy, the inverse matrix for the components of
such a form possesses different symmetry properties depending on
its parity,  corresponding, up to certain sign factors, to the
components of a bivector for an even form and to the coefficients
of a quadratic function
for an odd form.}. Given an even bivector field $P\in \fun(\Pi T^*M)$, a Poisson bracket
$\{f ,g \}_P$ of functions $f$ and $g$ is defined by the coordinate-free formula
\begin{equation}\label{newpoisson}
  \{f,g\}_P:=\{f,\{P,g\}\}
\end{equation}
where the brackets at the r.h.s. are the canonical Schouten brackets of multivector fields. The
Jacobi identity for $\{\ ,\ \}_P$ is equivalent to the identity $\{P,P\}=0$. Likewise, given an
odd quadratic Hamiltonian $S\in \fun(T^*M)$, it defines a Schouten structure  $\{\ ,\ \}_S$ by the
similar formula
\begin{equation}\label{newschouten}
  \{f,g\}_S:=\{f,\{S,g\}\}
\end{equation}
where at the r.h.s. there are the canonical Poisson brackets, and
the Jacobi identity for $\{\ ,\ \}_S$ is equivalent to the
identity $\{S,S\}=0$. The interrelation between Poisson and
Schouten structures is  summarized  in the table below (see also
Appendix B to~\cite{tv:as}). Notice that
$\{f,\{P,g\}\}=\{\{f,P\},g\}$, and $\{f,\{S,g\}\}=\{\{f,S\},g\}$.

\medskip
\begin{center}
{\renewcommand{\arraystretch}{1.3}
\begin{tabular}{|c|c|} \hline
  \rule{0pt}{18pt} A Poisson manifold $M$ & A Schouten manifold $M$  \\
  \hline
  \multicolumn{2}{|c|}{\rule{0pt}{21pt}$\displaystyle
  \{f,g\}_{P,S}= (-1)^{\at(\ft+1)}\der{f}{x^a}\{x^a,x^b\}_{P,S}\der{g}{x^b}$}\\
  \multicolumn{2}{|c|}{\small(same coordinate formula for even and odd bracket)}\\
  \hline
  $\{x^a,x^b\}_P=(-1)^{\at}\,P^{ab}$,  & $\{x^a,x^b\}_S=-(-1)^{\at}\,S^{ab}$,   \\
$P^{ab}$ the Poisson tensor (bivector field) & $S^{ab}$ the Schouten tensor (Hamiltonian)\\
$P^{ab}=(-1)^{(\at+1)(\bt+1)}P^{ba}$ & $S^{ab}=(-1)^{\at\bt}S^{ba}$\\
$\widetilde{P^{ab}}=\at+\bt$ & $\widetilde{S^{ab}}=\at+\bt+1$\\
$\displaystyle
P=\frac{1}{2}P^{ab}(x)x^*_bx^*_a\in \fun(\Pi T^*M)$ & $\displaystyle
S=\frac{1}{2}S^{ab}(x)p_bp_a\in \fun(T^*M)$\\
Jacobi for $\{\ ,\ \}_P$ $\Leftrightarrow$ $\{P,P\}=0$ & Jacobi for $\{\ ,\ \}_S$
$\Leftrightarrow$ $\{S,S\}=0$\\
Explicit formula: & Explicit formula:\\
$\displaystyle \{f,g\}_P=\{f,\{P,g\}\}$ & $\displaystyle
\{f,g\}_S=\{f,\{S,g\}\}$\\
\rule{0pt}{22pt}\parbox{6.5cm}{ \small Here $\{\ ,\ \}$ is the canonical  Schouten bracket on $\Pi %
T^*M$. $P$ is even.}  & \rule{0pt}{22pt}\parbox{6.5cm}{\small Here $\{\ ,\ \}$ is the canonical
Poisson bracket on $T^*M$.
\small  $S$ is odd.} \\
\rule{0pt}{-20pt} &  \rule{0pt}{-20pt}\\
   \hline
\end{tabular}
}
\end{center}

\medskip

(Notation: $x^a$ stand for local coordinates on $M$, $p_a$ for the
induced coordinates in the cotangent space, $x^*_a$ for the
induced coordinates in the anticotangent space. Here $\tilde
x^a=\at=\tilde p_a$, $\tilde x^*_a=\at+1$. The transformation law
for $p_a$ and  $x^*_a$ is the same.)

We shall later need the coordinate formulas for the canonical brackets.

On $T^*M$ the canonical even symplectic form is
$dp_adx^a=d(dx^ap_a)=d(p_adx^a)$.  The corresponding canonical
Poisson bracket:
\begin{equation}\label{canpoisson}
  \{f,g\}=(-1)^{\at(\ft+1)}\der{f}{p_a}\der{g}{x^a} - (-1)^{\at\ft}\der{f}{x^a}\der{g}{p_a},
\end{equation}
where $f,g\in \fun(T^*M)$, i.e.,  Hamiltonians on $M$. In particular,
$\{p_a,x^b\}=\delta_a^b=-(-1)^{\at}\{x^b,p_a\}$.

On $\Pi T^*M$ the canonical odd symplectic form  is $(-1)^{\at+1}dx^*_adx^a=d(dx^a
x^*_a)=d((-1)^{\at+1}x^*_adx^a)$. The corresponding canonical Schouten bracket:
\begin{equation}\label{canschouten}
  \{f,g\}=(-1)^{(\at+1)(\ft+1)}\der{f}{x^*_a}\der{g}{x^a} - (-1)^{\at(\ft+1)}\der{f}{x^a}\der{g}{x^*_a},
\end{equation}
where $f,g\in \fun(\Pi T^*M)$ (multivector fields on $M$). In particular,
$\{x^*_a,x^b\}=\delta_a^b=-\{x^b,x^*_a\}$.

The properties of the brackets~\eqref{newpoisson}
and~\eqref{newschouten}  are derived from the properties of the
canonical brackets. Crucial is that the functions on $M$
pulled-back to $T^*M$ or $\Pi T^*M$ have vanishing canonical
brackets. In particular, this implies antisymmetry of a new
bracket. In such description it is clear that this construction
can be put into an abstract setting and generalized (in
particular, iterated). It goes as follows.

\begin{de}[see~\cite{yvette:derived},\cite{yvette:loday}]
A \textit{Loday  algebra} of parity $s=0,1$, is a vector space $L$
endowed with a bilinear operation of parity $s$ (a \textit{Loday
bracket}) satisfying the condition
\begin{equation}\label{loday}
  [a,[b,c]]   =[[a,b]c]+(-1)^{(\at+s)(\bt+s)}[b,[a,c]].
\end{equation}
\end{de}
In other words, a Loday algebra is a crippled Lie superalgebra without the antisymmetry condition
and where the Jacobi identity remains in the form of a one-sided ``Leibniz identity'' for the
bracket. (Loday~\cite{loday:leibniz} introduced these algebras under the name of ``Leibniz
algebras''.)

Consider an odd derivation $D$ of a Loday algebra $L$. Introduce a new operation
$[a,b]_D:=(-1)^{\at}[Da,b]$. Obviously, it is a  bilinear operation of parity $s+1$.
\begin{thm}[\cite{yvette:derived}] If $D^2=0$, then the bracket $[\ ,\ ]_D$ defines a
new structure of a Loday algebra on the space $L$, of parity $s+1$.
\end{thm}
\begin{proof}
Consider
\begin{multline*}
[a,[b,c]_D]_D=(-1)^{\at+\bt}[Da,[Db,c]]=\\(-1)^{\at+\bt}\left(
[[Da,Db],c]+(-1)^{(\at+s+1)(\bt+s+1)}[Db,[Da,c]] \right)=\\
(-1)^{\at+\bt}\left( (-1)^{\at+s+1}[D[Da, b],c]+(-1)^{(\at+s+1)(\bt+s+1)}[Db,[Da,c]] \right),
\end{multline*}
because $D[Da,b]=[D^2a,b]+(-1)^{\at+s+1}[Da,Db]=-1)^{\at+s+1}[Da,Db]$, due to $D^2=0$. Thus it
equals (we continue)
\begin{multline*}
(-1)^{\at+\bt+s+1}[D(-1)^{\at}[Da,
b],c]+(-1)^{(\at+s+1)(\bt+s+1)}(-1)^{\bt}[Db,(-1)^{\at}[Da,c]]=\\
[[a,b]_D,c]_D+(-1)^{(\at+s+1)(\bt+s+1)}[b,[a,c]_D]_D.
\end{multline*}
\end{proof}

The new bracket (of the opposite parity) $[\ ,\ ]_D$ is called a
\textit{derived bracket}. Suppose $L$ is a Lie algebra  of parity
$s$ (cf. in~\cite{tv:poiss}), i.e., there is the antisymmetry
condition $[a,b]=-(-1)^{(\at+s)(\bt+s)}[b,a]$ in addition
to~\eqref{loday}. Without loss of generality set $s=0$. Then
$[a,b]_D+(-1)^{(\at+ 1)(\bt+
1)}[b,a]_D=(-1)^{\at}[Da,b]+(-1)^{\at\bt+\at+1}[Db,a]=(-1)^{\at}[Da,b]+[a,Db]=(-1)^{\at}D[a,b]$.
Thus on a subspace where $D[a,b]=0$ (or on a certain quotient) the
derived bracket gives a genuine Lie bracket of parity $s+1$.
(See~\cite{yvette:derived} for details and generalizations.
Without Loday algebras, derived brackets were also considered
in~\cite{tv:lectures93, tv:unpubl}.)

\begin{rem}\label{derivedandcartan}
The derived bracket construction is seen also in the formulas of
the Cartan calculus, where $[d,i_u]=L_u$ and
$[L_u,i_v]=(-1)^{\ut}i_{[u,v]}$ imply $i_{[u,v]}=[i_u,[d,i_v]]$,
which is a coordinate-free relation of the structure equation for
the commutator of vector fields with that for $d$.
\end{rem}

The construction of a Poisson/Schouten bracket from the canonical
brackets fits into this abstract model, with $D=\ad P$ or $D=\ad
S$, respectively. If we consider an arbitrary even multivector
field $P$ satisfying $\{P,P\}=0$, rather than a bivector field,
then the subalgebra $\fun(M)\subset\fun(\Pi T^*M)$ ceases to be
closed under the new bracket. Instead, we arrive at a
``Loday--Poisson'' structure on the whole of $\fun(\Pi T^*M)$.
Same holds for the Schouten case, if the restriction that the
Hamiltonian $S$ is quadratic in $p_a$ is dropped.

\begin{rem}
In the algebra of functions on a Schouten manifold $M$ endowed with a volume form there
is a natural ``odd Laplacian''
\begin{equation}\label{hoviksdelta}
  \Delta f:= \frac{1}{2}\,\div X_f,
\end{equation}
which makes it into a Batalin--Vilkovisky algebra. The
formula~\eqref{hoviksdelta} is due to
Khudaverdian~\cite{hov:delta}.  It provides an invariant geometric
setup for the $\Delta$-operator of the BV--formalism (see
also~\cite{hov:lap},\cite{tv:laplace1}). Notice an analogy with
the modular field on Poisson manifolds
(see~\cite{weinstein:modular}). For $\Pi T^*M$ considered with the
volume form $\boldsymbol\rho^2$, where $\boldsymbol\rho$ is a
volume form on $M$, the operator~\eqref{hoviksdelta} coincides
with the usual divergence of multivector fields $\delta$.
\end{rem}

\section{``Bi-'' structures and the Drinfeld double in a supermanifold framework}\label{bistructures}

Recall the following definition. For simplicity let us first consider the purely even situation.
Later we will give a different approach that  automatically handles the general case.

\begin{de}[Drinfeld]
A \textit{Lie bialgebra} is a Lie algebra $\mathfrak g$ such that
the dual space $\mathfrak g^*$ is also a Lie algebra and the two
structures are compatible in the sense that
\begin{equation}\label{bialgebra}
  \delta{\co}  \mathfrak g\to \wedge^2\mathfrak g
\end{equation}
is a $1$-cocycle of $\mathfrak g$ w.r.t. $\wedge^2 \ad$, where $\delta^*{\co}
\wedge^2\mathfrak g^*\to \mathfrak g^*$ is the bracket on $\mathfrak g^*$.
\end{de}

If $ e_i$ is a basis of $\mathfrak g$, $e^i$ the dual basis of $\mathfrak g^*$, then the
definition of a bialgebra is equivalent to the following relation between the structure constants
of $\mathfrak g$ and $\mathfrak g^*$ (the formula below is for a purely even $\mathfrak g$):
\begin{equation}\label{coorbialgebra}
  c_{jk}^i b_i^{nm}-c_{ji}^n b_k^{im}+c_{ji}^m b_k^{in}
  +c_{ki}^n b_j^{im}-c_{ki}^m b_j^{in}=0.
\end{equation}

\begin{thm}[Drinfeld]\label{drinfeld}
Let $\mathfrak g$ be a Lie bialgebra. Then on $\mathfrak g\oplus \mathfrak g^*$ there is a unique
structure of a Lie bialgebra  described by the properties:
\begin{enumerate}
  \item $\mathfrak g$ and $\mathfrak g^*$ are Lie subalgebras,
  \item The natural inner product in $\mathfrak g\oplus \mathfrak g^*$ is $\ad$-invariant,
  \item The cobracket $\delta$ for $\mathfrak g\oplus \mathfrak g^*$
  is given by the coboundary of the element
  $r=e_i\wedge e^i\in \wedge^2(\mathfrak g\oplus \mathfrak g^*)$.
\end{enumerate}
\end{thm}

A partial converse is given by the following theorem.
\begin{thm}[Manin] \label{manin}
Let $\mathfrak a, \mathfrak b$ be (finite-dimensional) Lie algebras such that the vector space
$\mathfrak d:=\mathfrak a \oplus \mathfrak b$ has a structure of a Lie algebra with an invariant
inner product. Suppose that
\begin{enumerate}
  \item $\mathfrak a$ and $\mathfrak b$ are isotropic subspaces in $\mathfrak d$,
  \item $\mathfrak a$ and $\mathfrak b$ as Lie  algebras are subalgebras in $\mathfrak d$.
\end{enumerate}
Then $\mathfrak b \cong \mathfrak a^*$, and $\mathfrak a$, $\mathfrak b\cong \mathfrak a^*$ are
dual Lie bialgebras. The Lie algebra structure on $\mathfrak d$ coincides with the one given by
Drinfeld's theorem.
\end{thm}

The space $\mathfrak{d}=\mathfrak{d(g)}:=\mathfrak g\oplus \mathfrak g^*$ with the bialgebra
structure is called the (classical) \textit{Drinfeld double} of the bialgebra $\mathfrak g$; the
triple  of Lie algebras $(\mathfrak a, \mathfrak b, \mathfrak d=\mathfrak a \oplus \mathfrak b)$ is
called a \textit{Manin triple}. Explicit formulas for the bracket in $\mathfrak d$:
\begin{align}
  [e_i,e_j]&=c_{ij}^k e_k, \label{doublealg1}\\
  [e^i,e^j]&=b^{ij}_k e^k,\\
  [e_i,e^j]&=b^{jk}_i e_k - c_{ik}^j e^k.\label{doublealg3}
\end{align}
The cobracket $\delta{\co}  \mathfrak d\to \wedge^2\mathfrak d$ is given by
\begin{align}\label{doublecoalg1}
  \delta(e_i) & = b_i^{jk}e_j\wedge e_k\\
  \delta(e^i) & = -c^i_{jk}e^j\wedge e^k. \label{doublecoalg2}
\end{align}
(Due to the natural inner product, the cobracket~(\ref{doublecoalg1},\ref{doublecoalg2}) gives rise
to the second bracket in $\mathfrak d$: $[e_i,e_j]_{(2)}=-c_{ij}^k e_k$, $[e^i,e^j]_{(2)}=b^{ij}_k
e^k$, $[e_i,e^j]_{(2)}=0$.)

\bigskip
Now, to get a different formulation suitable for further
generalizations, we first have to describe a Lie algebra structure
geometrically. To this end, let $\mathfrak g=\mathfrak g_0\oplus
\mathfrak g_1$ be a vector space (from now on we return to the
general super case). Consider it together with its neighbors. (The
term ``neighbor'' for super vector spaces is due to Manin.) We
have a tetrahedron\footnote{The picture of a tetrahedron was
suggested to me by Alan~Weinstein instead of the three-pointed
star with center at $\mathfrak g$ that I used initially.}:

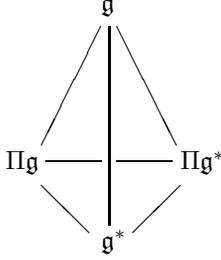
\begin{center}
{ \unitlength=3pt
\begin{picture}(0,35)
\put(0,35)
    {\begin{picture}(0,0) \put(0,0){$\mathfrak g$} \put(-12,-20){$\Pi \mathfrak
    g$} \put(10,-20){$\Pi \mathfrak g^*$} \put(0,-30){$\mathfrak g^*$} {\thicklines
    \put(1,-2){\line(0,-1){25}} }\put(-7,-19){\line(1,0){7}} \put(2,-19){\line(1,0){7}}
    \put(0,-2){\line(-1,-2){7.5}}\put(2,-2){\line(1,-2){7.5}}
    \put(-7.5,-22){\line(1,-1){6}}\put(10,-22){\line(-1,-1){6}}
    \end{picture}}
\end{picture}
}
\end{center}

\begin{thm}\label{liealgebra}
A Lie (super)algebra structure in the vector space $\mathfrak g$ is equivalent to either of:
\begin{enumerate}
  \item A linear Poisson bracket on the manifold  $\mathfrak g^*$ (Lie coalgebra).
  \item A linear Schouten bracket on the manifold  $\Pi \mathfrak g^*$ (Lie anticoalgebra).
  \item A quadratic homological vector field on the manifold $\Pi \mathfrak g$ (Lie antialgebra).
\end{enumerate}
\end{thm}

Recall that a \textit{homological vector field} on a supermanifold
$M$ is an odd vector field (=odd derivation of the algebra of
functions) $\hat Q$ such that $[\hat Q,\hat Q]=2\hat Q^2=0$. A
basis $e_i$ in $\mathfrak g$ gives rise to linear coordinates on
the supermanifolds $\Pi \mathfrak g$, $\mathfrak g^*$, and $\Pi
\mathfrak g^*$, which we denote $\xi^i$, $x_i$  and $\xi_i$,
respectively. If $\tilde e_i=:\itt$, then $\tilde \xi^i=\itt+1$,
$\tilde x_i=\itt$, $\tilde \xi_i=\itt+1$.

\begin{ex}
The homological field corresponding to a Lie superalgebra $\mathfrak g$ with the structure
constants $c_{ij}^k$ (i.e., $[e_i,e_j]=c_{ij}^k e_k$) is
\begin{equation}\label{lie}
  \hat Q=\frac{1}{2}(-1)^{\jt}\xi^j\xi^ic_{ij}^k\,\der{}{\xi^k}.
\end{equation}
The condition $\hat Q^2=0$ is equivalent to the Jacobi identity
for $c_{ij}^k$. An element $X\in \mathfrak g$ (here we treat
$\mathfrak g$ as a vector space, not as a manifold) corresponds to
a vector field on $\Pi \mathfrak g$ (of parity opposite to that of
$X$)  with constant coefficients: $X=X^ie_i\mapsto
i_X:=(-1)^{\Xt}X^i\lder{}{\xi^i}$. The odd linear map $i{\co}
\mathfrak g\to \Vect(\Pi \mathfrak g)$ is,  obviously, injective.
The Lie bracket in $\mathfrak g$ can be reconstructed from the
commutator of vector fields on $\Pi \mathfrak g$ by the ``derived
bracket'' formula:
\begin{equation}\label{liebracket}
  i_{[X,Y]}=[i_X,[\hat Q,i_Y]]=[[i_X,\hat Q],i_Y].
\end{equation}
Here at the l.h.s. $[\ ,\ ]$ stands for the Lie bracket in $\mathfrak g$ with the
structure constants $c_{ij}^k$, and at the r.h.s. $[\ ,\ ]$ stands for the canonical Lie
bracket (commutator) of vector fields. (Notice that vector fields with constant
coefficients have zero commutator.) One can recognize in~\eqref{liebracket} an analog of
Cartan's formulas, see Remark~\ref{derivedandcartan}.

\begin{ex} For  the general linear algebra $\mathfrak{gl}(n)$ we have the field
\begin{equation}
  \hat Q=-\sum \x^{ik}\x^{kj}\der{}{\x^{ij}}\,.
\end{equation}
It is tangent to linear submanifolds that correspond to Lie subalgebras in
$\mathfrak{gl}(n)$. Same formula works in the super case.
\end{ex}

In a more conventional language, the vector field $\hat Q$ on $\Pi
\mathfrak g$ coincides with the Chevalley--Eilenberg differential
in $C^*(\mathfrak g)=C^*(\mathfrak g;\RR)$, for ordinary Lie
algebras. For Lie superalgebras this  can be taken as a convenient
definition.


\end{ex}
\begin{ex}Consider a general odd (formal) vector field on $\mathbb R^{n|m}$:
\begin{equation}\label{holie}
  \hat Q=\left(Q_0^k+\xi^iQ^k_i+\frac{1}{2}\xi^j\xi^iQ_{ij}^k+
  \frac{1}{3!}\xi^l\xi^j\xi^iQ_{ijl}^k+\ldots\right)\der{}{\xi^k}.
\end{equation}
The coefficients $Q_0^k$, $Q^k_i$, $Q_{ij}^k$, $Q_{ijl}^k$, \dots
define a sequence of $N$-ary operations ($N=0,1,2,3,\dots$) on the
vector space $\mathbb R^{m|n}$, and the condition $\hat Q^2=0$
expands to a linked sequence of ``Jacobi identities''. If only the
quadratic term in the Maclaurin expansion~\eqref{holie} is
present, we return to the case of a Lie algebra. The general case
is a  \textit{strong homotopy Lie algebra}
(\textit{$L_{\infty}$-algebra}) (due to Stasheff,
see~\cite{stasheff:shla93}).
\end{ex}
\begin{ex}
If instead of commutative algebras of functions on $\R{n|m}$ one
considers free noncommutative associative algebras and their odd
derivations of square zero (analog of homological fields), then  a
homotopy version of associative algebras (so called
\textit{$A_{\infty}$-algebras}) will be
obtained~\cite{stasheff:intrinsic}. Again, the purely quadratic
case reproduces ordinary (associative) algebras.
\end{ex}

Now we can say  what is, in this language, a general ``bi-'' structure for a Lie algebra. Consider
the above tetrahedron of vector spaces. A ``bi-'' structure for $\mathfrak g$ is exactly a Lie
algebra structure for two of its vertices, with a certain compatibility condition. In view of
Theorem~\ref{liealgebra}, up to renaming, the list of the corresponding geometric structures is
exhausted by three cases:

\begin{center}
{\renewcommand{\arraystretch}{1.3}
\begin{tabular}{|r|l|l|l|}
\hline
  \rule{0pt}{18pt} I. & $\mathfrak g$, $\mathfrak g^*$ Lie algebras
  & \rule{0pt}{22pt}\parbox{6.0cm}{On $\Pi
  \mathfrak g$ (or $\Pi \mathfrak g^*$) there are homological field $\hat Q$
  and odd bracket.} & $QS$-manifold\\
  \hline
    II. & $\mathfrak g$, $\Pi \mathfrak g^*$ Lie algebras
  & \rule{0pt}{22pt}\parbox{6.0cm}{On $\mathfrak g^*$ (or $\Pi \mathfrak g$)
  there are homological field $\hat Q$ and even bracket.} & $QP$-manifold\\
  \hline
  III. & $\mathfrak g$, $\Pi \mathfrak g$ Lie algebras
  & \rule{0pt}{22pt}\parbox{6.0cm}{On $\mathfrak g^*$
  (or $\Pi \mathfrak g^*$) there are even and odd brackets.\vphantom{Q}} & $PS$-manifold\\
  \hline
\end{tabular}
}
\end{center}
One should explain the compatibility conditions. For a pair consisting of a Schouten or Poisson
bracket and a homological vector field there is only one reasonable condition, that the field
$\hat Q$ should be a derivation of the bracket. We call Poisson (Schouten) manifolds endowed with a
homological field which is a derivation of the bracket \textit{$QP$-manifolds}
(\textit{$QS$-manifolds,} respectively).

\begin{thm}\label{qsstructure}
A $QS$-manifold structure on $\Pi \mathfrak g$ with a linear
Schouten bracket and a quadratic field $\hat Q$ is equivalent to a
Lie bialgebra structure for $\mathfrak g$.
\end{thm}
\begin{proof}
As noted above, the vector field $\hat Q$ corresponds to the
Chevalley--Eilenberg differential for the trivial  representation.
The Lie derivative $L_{\hat Q}$ corresponds to the
Chevalley--Eilenberg differential for other (nontrivial)
representations. Thus the condition that the cobracket on
$\mathfrak g$ is a cocycle is equivalent to $L_{\hat Q}S=0$, for
the Hamiltonian $S$ that defines the corresponding Lie--Schouten
structure on $\Pi \mathfrak g$. On the other hand, we have, by a
straightforward calculation,
\begin{equation}
  {\hat Q}\{f,g\}_S=\{\hat Q f,g\}_S+(-1)^{\ft+1}\{f,\hat Q g\}_S+
  (-1)^{\ft}\{f,\{L_{\hat Q}S,g\}\},
\end{equation}
for the bracket $\{\ ,\ \}_S$ defined by $S$. Hence, if $L_{\hat
Q}S=0$, then the derivation property holds. Conversely, if the
derivation property holds, then for all $f$, $g$ $\{f,\{L_{\hat
Q}S,g\}\}=0$, and this implies $L_{\hat Q}S=0$ by the virtue of
the non-degeneracy of the canonical Poisson bracket.
\end{proof}

In particular, this is the way to get the correct formulas like~\eqref{coorbialgebra} in the super
case $\mathfrak g=\mathfrak g_0\oplus \mathfrak g_1$. Practically, this theorem can be used
instead of the definition in this case.

What about the $QP$-manifolds? They correspond to an ``odd version'' of the Lie bialgebra notion.

\begin{de}
An \textit{odd Lie bialgebra} is a Lie superalgebra $\mathfrak g$ such that the space $\Pi
\mathfrak g^*$ is also a Lie superalgebra so that on  the antialgebra $\Pi \mathfrak g$ is induced
a $QP$-structure (with linear Poisson bracket and quadratic homological field).
\end{de}

We shall elaborate this in Section~\ref{odd}.

\begin{rem}
It should be emphasized that while a Lie algebra $\mathfrak g$ can be viewed via three
geometric manifestations, listed in Theorem~\ref{liealgebra}, for a \mbox{``bi-''} case
there is just one geometric picture. A $QS$-structure with linear bracket and quadratic
field (realized, equivalently, either on $\Pi \mathfrak g$ or on $\Pi\mathfrak g^*$)
corresponds to a Drinfeld's Lie bialgebra, and a $QP$-structure with linear bracket and
quadratic field (realized, equivalently, either on $\Pi \mathfrak g$ or on $\mathfrak
g^*$) corresponds to an odd bialgebra.
\end{rem}

{\footnotesize  We reserve the name \textit{$PS$-manifold} for a manifold with a pair of
brackets of opposite parity. The precise compatibility condition in this case is not
obvious. We refrain from discussing it here. Geometry of spaces with even and odd bracket
was studied, in particular, in~\cite{hov:delta}. For  $\mathfrak g\oplus\Pi \mathfrak g$,
an odd endomorphism replaces  the pairing.

}

\smallskip
Now, we need the description of Drinfeld's double in this
language. Start from the space $\Pi \mathfrak g$. Notice that $\Pi
(\mathfrak g\oplus \mathfrak g^*)\cong T^*\Pi \mathfrak g$. The
natural (symmetric) inner product on $\mathfrak g\oplus \mathfrak
g^*$ corresponds to the canonical symplectic structure on $T^*\Pi
\mathfrak g$. Denote coordinates on $\Pi \mathfrak g$ by $\xi^i$
and their conjugate momenta by $\xi_j$. As we know from
Section~\ref{psalgebras}, a linear Schouten bracket on $\Pi
\mathfrak g$ is specified by an odd Hamiltonian of the form
$S=(1/2)\xi^lS_l^{nm}\xi_m\xi_n\in\fun(T^*\Pi \mathfrak g)$. On
the other hand, a vector field $\hat
Q=(1/2)\xi^j\xi^iQ_{ij}^k\lder{}{\xi^k}\in\Vect(\Pi \mathfrak g)$
induces the Hamiltonian $Q=(1/2)\xi^j\xi^iQ_{ij}^k\x_k\in
\fun(T^*\Pi \mathfrak g)$.

\begin{thm}\label{thmbialgebradouble}
The following statements hold for $T^*\Pi \mathfrak g\cong \Pi (\mathfrak g\oplus \mathfrak g^*)$:

{\em 1.}  The $QS$-condition for $\Pi \mathfrak g$ with the field $\hat Q$ and the
Schouten tensor $S$ is equivalent to the vanishing of the Poisson bracket
\begin{equation}\label{invbialgebra}
        \{Q,S\}=\frac{1}{4}\{\xi^j\xi^iQ_{ij}^k\x_k,\xi^lS_l^{nm}\xi_m\xi_n\}=0.
\end{equation}

{\em 2.}  The function $=Q+S$ on $T^*\Pi \mathfrak g$ satisfies $\{Q+S,Q+S\}=0$, so the
Hamiltonian vector field $\hat Q_D:=X_{Q+S}\in \Vect(T^*\Pi \mathfrak g)$  is
  homological. It corresponds to the Lie algebra structure in Drinfeld's double
  $\mathfrak d=\mathfrak g\oplus \mathfrak g^*$.
\end{thm}
\begin{proof}
To prove part 1, notice that $\{Q,S\}=L_{\hat Q}S$, thus the equation $\{Q,S\}=0$ is equivalent to
the $QS$-condition (see the proof of Theorem~\ref{qsstructure}). To prove part 2, notice that the
equation $\{Q+S,Q+S\}=0$ implies the equation $\hat Q_D^2=X_{Q+S}^2=0$ for the corresponding
Hamiltonian field. Finally, to show that $\hat Q_D$ gives  exactly the Lie bracket in Drinfeld's
double, we compare the explicit formulas. For $X_Q$ we obtain, using the general
formula~\eqref{canpoisson}:
\begin{multline*}
    X_Q=\{Q,\ \}=
    (-1)^{(\kt+1)(\Qt+1)}\der{Q}{\xi_k}\der{}{\xi^k}-(-1)^{(\kt+1)\Qt}\der{Q}{\xi^k}\der{}{\xi_k}\\
    = \der{Q}{\xi_k}\der{}{\xi^k}+(-1)^{\kt}\der{Q}{\xi^k}\der{}{\xi_k}=
    \frac{1}{2}\x^j\x^iQ_{ij}^k\der{}{\x^k}+(-1)^{\kt}\x^iQ_{ik}^j\x_j\der{}{\x_k}
\end{multline*}
Similarly, we obtain
\begin{multline*}
    X_S=\{S,\ \}
    = \der{S}{\xi_k}\der{}{\xi^k}+(-1)^{\kt}\der{S}{\xi^k}\der{}{\xi_k}\\=
    \x_n\x^lS^{nk}_l\der{}{\x^k}+(-1)^{\kt}\frac{1}{2}S^{nm}_k\x_m\x_n\der{}{\x_k}
\end{multline*}
Thus $\hat Q_D=X_Q+X_S$ reads
\begin{equation}\label{qd}
  \hat Q_D=\left(\frac{1}{2}\x^j\x^iQ_{ij}^k+\x_j\x^iS^{jk}_i\right)\der{}{\x^k}
  +
  (-1)^{\kt}\left(\frac{1}{2}S^{ij}_k\x_j\x_i+\x^iQ_{ik}^j\x_j\right)\der{}{\x_k},
\end{equation}
and in the case of a purely even $\mathfrak g$ it is reduced to
\begin{equation*}
  \left(\frac{1}{2}\x^j\x^iQ_{ij}^k+\x_j\x^iS^{jk}_i\right)\der{}{\x^k}
  +
  \left(\frac{1}{2}S^{ij}_k\x_j\x_i-\x_j\x^iQ_{ik}^j\right)\der{}{\x_k},
\end{equation*}
which reproduces formulas~(\ref{doublealg1}--\ref{doublealg3}) for
the Lie bracket.
\end{proof}

\begin{rem}
It should be noted that the description of the double of a Lie
bialgebra in terms of a pair of commuting Hamiltonians was first
done in~\cite{lecomte:roger}. A similar description  for Lie
``quasi-bialgebras'' was obtained in~\cite{yvette:jacobian}.
\end{rem}

Theorem~\ref{thmbialgebradouble} gives for $T^*\Pi \mathfrak g$ a
$Q$-structure corresponding to the Lie bracket in the double. To
get on $T^*\Pi \mathfrak g$ an $S$-structure corresponding to the
cobracket, consider the cotangent bundle of $T^*\Pi \mathfrak g$.
Denote by $\pi_i,\pi^i$ the conjugate momenta for $\xi^i,\xi_i$,
respectively. Notice that $\pi_i$ transforms as $\x_i$, and
$\pi^i$ as $(-1)^{\itt}\xi^i$ (see a more general statement in
Section~\ref{bialgebroids}). Hence
$r=-(-1)^{\itt}\pi^i\pi_i=\pi_i\pi^i$ is an invariantly defined
function on $T^*T^*\Pi \mathfrak g$.
Consider the following linear map $\fun(T^*\Pi\mathfrak g)\to\fun(T^*T^*\Pi\mathfrak g)$,
$f\mapsto \bar f:=\{p_{X_f},r\}$.

\begin{lm}
\label{lfbar} For an arbitrary function $f\in \fun(T^*\Pi\mathfrak g)$ there is an explicit formula
\begin{equation}\label{fbar}
  \bar f=\{p_{X_f},r\}=(-1)^{(\itt+\jt)\ft}\dder{f}{\x^j}{\x^i}\p^i\p^j
  -(-1)^{(\itt+\jt)(\ft+1)}\dder{f}{\x_j}{\x_i}\p_i\p_j.
\end{equation}
\end{lm}
\begin{proof}
Notice, first, that the bracket of an arbitrary $g\in \fun(T^*T^*\Pi\mathfrak g)$ with
$r=\p_i\p^i$  equals
\begin{equation}\label{skobkasr}
  \{g,r\}=-(-1)^{(\jt+1)\gtt}\der{g}{\x^j}\p^j-(-1)^{(\jt+1)(\gtt+1)}\der{g}{\x_j}\p_j.
\end{equation}
Now,  for an arbitrary $f\in \fun(T^*\Pi\mathfrak g)$ we get the function
\begin{equation}\label{pxf}
  p_{X_f}=(-1)^{(\itt+1)(\ft+1)}\der{f}{\x_i}\p_i-(-1)^{(\itt+1)\ft}\der{f}{\x^i}\p^i
\end{equation}
on $T^*T^*\Pi\mathfrak g$. Substituting~\eqref{pxf}
into~\eqref{skobkasr}, we obtain formula~\eqref{fbar} by a direct
simplification (which we omit). Notice that the terms containing
``cross'' products like $\pi_i\p^j$  remarkably cancel.
\end{proof}

\begin{thm}\label{koskobkavduble}
Consider  $S_D:=(1/2)\{Q_D,r\}\in \fun(T^*T^*\Pi \mathfrak g)$, where $Q_D
\in
\fun(T^*T^*\Pi \mathfrak g)$  corresponds to the vector field $\hat Q_D$. Then $S_D$ is a Schouten
tensor for $T^*\Pi \mathfrak g$,  and together with  the field $\hat Q_D$ it gives a
$QS$-structure on $T^*\Pi \mathfrak g$. The Schouten bracket  specified by $S_D$ corresponds  to
the cobracket in Drinfeld's double  $\mathfrak d=\mathfrak g\oplus \mathfrak g^*$.
\end{thm}
\begin{proof}
We have $S_D=(1/2)(\widebar{Q+S})=(1/2)(\bar Q+\bar S)$, in the notation of Lemma~\ref{lfbar}.
Applying the lemma to the functions $Q$ and $S$, we obtain the formulas
\begin{align}
  \bar Q&=(-1)^{\itt+\jt}Q_{ij}^k\x_k\p^j\p^i=(-1)^{\jt}\p^j(-1)^{\itt}\p^iQ_{ij}^k\x_k
  \label{barq}\\
  \bar S&=-\x^k S_k^{ij}\p_j\p_i.
\end{align}
Here in~\eqref{barq} we changed the order and separated the
factors $(-1)^{\itt}\p^i$. The point is that the correspondence
$\x^i\leftrightarrow (-1)^{\itt}\p^i$ respects the Poisson
bracket: $\{\x_i,\x^j\}_{T^*\Pi\mathfrak
g}=\delta_i^j=\{\x_i,(-1)^{\jt}\p^j\}_{T^*T^*\Pi\mathfrak g}$. In
the same way, $\{\x_i,\x^j\}_{T^*\Pi\mathfrak
g}=\delta_i^j=\{\p_i,\x^j\}_{T^*T^*\Pi\mathfrak g}$. Hence $\{\bar
Q,\bar Q\}=4\{Q,Q\}=0$, $\{\bar S,\bar S\}=4\{S,S\}=0$ by the
virtue of these correspondences, and $\{\bar Q,\bar S\}=0$ because
$\bar Q$ and $\bar S$ do not contain conjugate variables. Thus,
$\{S_D,S_D\}=0$. The identity $\{Q_D,S_D\}=0$ holds by the
construction of $S_D$. To notice that $S_D$ corresponds to the
cobracket in Drinfeld's double (up to a common sign), it suffices
to compare the explicit formula obtained
\begin{equation}\label{sdnabialgebre}
  S_D=\frac{1}{2}(-1)^{\itt+\jt}Q_{ij}^k\x_k\p^j\p^i-\frac{1}{2}\x^kS_k^{ij}\p_j\p_i
\end{equation}
with formulas~(\ref{doublecoalg1},\ref{doublecoalg2}) in the even
case.
\end{proof}

(It is clear now  how to get an ``odd analog'' of Drinfeld's  double. Let $\mathfrak g$ be an
{odd bialgebra}, so $\Pi \mathfrak g$ is a $QP$-manifold with the quadratic $\hat Q$ and
linear bracket. Instead of $T^*\Pi \mathfrak g$ we consider the supermanifold $\Pi T^*\Pi
\mathfrak g\cong \Pi \mathfrak g\oplus \mathfrak g^*$. In a way similar to the above, it can
be proved that $\Pi T^*\Pi \mathfrak g$ has a desired $QP$-structure, corresponding to an odd
bialgebra structure in $\mathfrak g\oplus \Pi\mathfrak g^*$. The space $\mathfrak g\oplus
\Pi\mathfrak g^*$ with this structure will be called the \textit{odd double} of an odd
bialgebra $\mathfrak g$. We will return to this construction in Section~\ref{odd} and give
details.)

\begin{rem}
Drinfeld's classical double is, by the construction, a 
\textit{coboundary Lie bialgebra}~(see \cite{charipressley}). That
means that $S_D=\{Q_D,r/2\}$, in our language. Suppose that for an
arbitrary Lie algebra $\mathfrak g$  an even Hamiltonian
$r\in\fun(T^*\Pi\mathfrak g)$ of the form $r=r^{ij}\p_i\p_j$ (with
constant $r^{ij}$) is considered. The Poisson bracket with the odd
Hamiltonian $Q=(1/2)\xi^j\x^iQ_{ij}^k\p_k$ corresponding to the
Lie algebra structure, $\{Q,r/2\}=:S$, is quadratic in $\p_i$:
$S=(1/2)\x^i S_i^{jk}\p_k\p_j$. The linear odd bracket on
$\Pi\mathfrak g$ defined as $\{f,g\}_S=\{f,\{S,g\}\}$ is
automatically compatible with $Q$. It satisfies the Jacobi
identity (thus making $\mathfrak g$ into a Lie bialgebra) if and
only if the Poisson bracket $\{\{Q,r\},\{Q,r\}\}$ vanishes. Since
$\{Q,\{r,\{Q,r\}\}\}=\{\{Q,r\},\{Q,r\}\}$, this is equivalent to
\begin{equation}\label{gyb}
  \{Q,\{r,r\}_Q\}=0
\end{equation}
where the inner bracket denotes the operation
\begin{equation}
  r_1,r_2\mapsto \{r_1,r_2\}_Q:=\{r_1,\{Q,r_2\}\}=\{\{r_1,Q\},r_2\}
\end{equation}
on Hamiltonians $r_1$, $r_2$ depending only on $\p_i$ (this is
exactly the Lie--Schouten bracket on $\Pi \mathfrak g^*$). The
equation~\eqref{gyb} is satisfied if, in particular,
\begin{equation}\label{yb}
  \{r,r\}_Q=0.
\end{equation}
This is the \textit{classical Yang--Baxter equation} for $r$,
while~\eqref{gyb} is the ``generalized'' classical Yang--Baxter
equation. Notice that $r=\p_i\p^i$ for Drinfeld's double
satisfies~\eqref{gyb}, but it does not satisfy~\eqref{yb},
$\{r,r\}_{Q_D}$ being a certain canonical $\hat Q_D$-invariant
cubic expression in $\p_i,\p^j$ (in the standard language $r$ is
said to be a solution of the ``modified'' classical Yang--Baxter
equation).
\end{rem}

\section{Doubles for Lie bialgebroids}\label{bialgebroids}
In this section we review the results of
D.~Roytenberg~\cite{roytenberg:thesis} and
A.~Vaintrob~\cite{vaintrob:algebroids}, who independently applied
super methods to Lie bialgebroids. As noted before, we give a
``superized'' version, i.e., consider from the beginning super Lie
(bi)algebroids over supermanifolds. Lie bialgebroids were
originally introduced by Mackenzie and Xu~\cite{mackenzie:bialg}.
Important results, which particularly nicely fit into the super
framework (but were originally formulated without it), are due to
Y.~Kosmann-Schwarzbach~\cite{yvette:exact}.

Consider a Lie algebroid $E\to M$ as defined in Section~\ref{psalgebras}. In direct analogy with
the tangent bundles and Lie algebras, it gives rise to a differential algebra, which is
$C^{\infty}(\Pi E)$ with the homological vector field $\hat Q=\hat Q_E\in\Vect(\Pi E)$,
\begin{equation}\label{algebroid}
  \hat Q  =\x^iQ_i^a(x)\,\der{}{x^a}+\frac{1}{2}\,\x^j\x^iQ_{ij}^k(x)\,\der{}{\x^k}.
\end{equation}
(Here we use the notation $x^a$ for coordinates in the base and
$\x^i$ for coordinates in the fiber of $\Pi E$.) In the purely
even case, $C^{\infty}(\Pi E)$ is the algebra of  smooth sections
of the exterior bundle $\fun(M,\Lambda(E^*))$. The differential
$d_E=\hat Q_E{\co} \fun(M,\Lambda(E^*))\to \fun(M,\Lambda(E^*))$
is analogous to the de Rham and Chevalley--Eilenberg differentials
(which are both due to Cartan). The coefficients $Q_{ij}^k(x)$
describe  the Lie brackets of a local basis of sections in $E$,
while $Q_i^a(x)$ correspond to the anchor map $a{\co}
\fun(M,E)\to \Vect(M)$. The bracket $[\ ,\ ]_E$ and the anchor are
recovered from $\hat Q$ by the formulas similar
to~\eqref{liebracket}:
\begin{align}
  a(u)f &= (-1)^{\ut}[i_u,[\hat Q,f]]= (-1)^{\ut}[[i_u,\hat Q],f]\\
  i_{[u,v]_E} & =[i_u,[\hat Q,i_v]]= [[i_u,\hat Q],i_v]
\end{align}
where $i_u=(-1)^{\ut}u^i(x)\lder{}{\x^i}$ for
$u=u^i(x)e_i\in\fun(M,E)$. In particular, for the basis sections
we have
\begin{equation}
  a(e_i)=Q_i^a(x)\,\der{}{x^a}, \quad [e_i,e_j]_E=(-1)^{\jt}Q_{ij}^k(x)e_k.
\end{equation}
The identity $\hat Q^2=0$ compactly contains all the properties of
the anchor and the bracket. Thus, we can take it as an equivalent
definition of a Lie algebroid structure. A vector bundle endowed
with $\hat Q$ of the form~\eqref{algebroid} such that $\hat Q^2=0$
is, by definition, a \textit{Lie antialgebroid}. $\Pi E\to M$ is
the Lie antialgebroid corresponding to the Lie algebroid $E\to M$.

The Lie algebroid structure in $E$ also defines a Schouten bracket
in the algebra $\fun(\Pi E^*)$, which in the purely even case is
the algebra of  smooth sections $\fun(M,\Lambda(E))$, and a
Poisson bracket in the algebra $\fun(E^*)$. Thus, as for Lie
algebras, there are a $Q$-manifold, an $S$-manifold and a
$P$-manifold associated with a Lie algebroid.

Return for a moment to a purely even situation. Let the dual bundle $E^*\to M$ also have
a structure of a Lie algebroid. Then, there is a differential
\begin{equation*}
    d_{E^*}{\co} \fun(M,\Lambda(E))\to \fun(M,\Lambda(E)).
\end{equation*}
It can be applied, in particular, to sections of $E\to M$
considered as sections of the exterior bundle.

\begin{de}[Mackenzie--Xu~\cite{mackenzie:bialg}]
A \textit{Lie bialgebroid} over a manifold $M$ is a Lie algebroid $E\to M$ for which the dual
bundle $E^*\to M$ is also a Lie algebroid so that the differential $d_{E^*}$ satisfies the
derivation property for the bracket of sections of $E$.
\end{de}

It was shown by Y.~Kosmann-Schwarzbach~\cite{yvette:exact} that it
is equivalent and more convenient  to  require that $d_{E^*}$ is a
derivation of the Schouten bracket in the whole algebra
$\fun(M,\Lambda(E))$. It can also be shown (though it is not
straightforward at this point) that $E$ and $E^*$ here can be
interchanged. Assuming this,
the definition can be put into a super language (allowing to cover the general case) as
follows:

\begin{de}[equivalent]\label{bialgebroid}
A (super) \textit{Lie bialgebroid} over a supermanifold $M$ is a vector bundle $E\to M$
such that the space of the opposite vector bundle $\Pi E\to M$ is a $QS$-manifold, with
the homological field of the form
\begin{equation}\label{bialgebroid1}
  \hat Q  =\x^iQ_i^a(x)\,\der{}{x^a}+\frac{1}{2}\,\x^j\x^iQ_{ij}^k(x)\,\der{}{\x^k}
\end{equation}
and the non-vanishing  Schouten brackets of coordinates of the form
\begin{equation}\label{bialgebroid2}
  \{\x^i,x^a\}_S= (-1)^{\itt}S^{ia}(x), \quad   \{\x^i,\x^j\}_S
  =(-1)^{\itt}\x^k S^{ij}_k(x).
\end{equation}
That means that the bracket is given by the Hamiltonian $S=S_E\in\fun(T^*\Pi E)$,
\begin{equation}\label{bialgebroid3}
  S=S^{ia}(x)p_a\p_i +\frac{1}{2}\x^k S^{ij}_k(x)\p_j\p_i.
\end{equation}
\end{de}

\begin{thm}\label{qsravno0}
$E$ is a Lie bialgebroid if and only if the Poisson bracket $\{Q,S\}$ vanishes, where
$Q=p(\hat Q)\in\fun(T^*\Pi E)$ is the Hamiltonian corresponding to the vector field $\hat
Q$.
\end{thm}
\begin{proof}
Similarly to Theorem~\ref{qsstructure}: the condition of a Lie bialgebroid is equivalent
to $L_{\hat Q_E}S_E=0$, and $L_{\hat Q_E}S_E=\{Q,S\}$.
\end{proof}

This definition 
is not manifestly symmetric with respect to $E$ and $E^*$. A
symmetric description can be achieved with the help of the
following theorem  (see Mackenzie and Xu~\cite{mackenzie:bialg};
we provide a proof in the Appendix).

\begin{thm} \label{testar}
For an arbitrary vector bundle $E\to M$ there is a natural diffeomorphism of the
cotangent bundles $F{\co} T^*E\cong T^*E^*$, which preserves the symplectic structure. In
coordinates:
\begin{equation*}
  F{\co}  (x^a,y^i,p_a,p_i)\mapsto (x^a,y_i,p_a,p^i)=(x^a,p_i,p_a,-(-1)^{\itt}y^i).
\end{equation*}
Here $y^i$ stand for (left) coordinates in the fiber of $E$, $y_i$
stand for the contragredient right coordinates in the fiber of
$E^*$, and $p_i$, $p^i$ for the respective conjugate momenta.
\end{thm}

Consider a Lie algebroid $\Pi E$. Let  the vector field $\hat
Q=\hat Q_E\in\Vect (\Pi E)$ be  given by~\eqref{bialgebroid1}. The
corresponding Hamiltonian $Q_E\in\fun(T^*\Pi E)$ is
\begin{equation}
  Q_E  =\x^iQ_i^a(x)\,p_a+\frac{1}{2}\,\x^j\x^iQ_{ij}^k(x)\,\p_k.
\end{equation}
We shall apply Theorem~\ref{testar} to the vector bundle $\Pi E$. Suppose $E^*$ is also a
Lie algebroid. Let $\hat Q_{E^*}\in\Vect (\Pi E^*)$ be the homological field defining the
algebroid structure in $E^*$ and the corresponding Hamiltonian $Q_{E^*}\in\fun(T^*\Pi
E^*)$ be
\begin{equation}
  Q_{E^*}  =\x_i Q^{ia}(x) p_a+\frac{1}{2}\x_j\x_iQ^{ij}_k(x)\p^k.
\end{equation}
Here $\x_i$ with lower indices denote the coordinates in the fiber
of $\Pi E^*$ and $\p^i$ the conjugate momenta. In particular, the
anchor and  bracket in $E^*$ will be given by
\begin{equation}\label{brackinestar}
  a(e^i)=Q^{ia}(x)\der{}{x^a}, \quad
  [e^i,e^j]_{E^*}=(-1)^{\jt}Q^{ij}_k(x)e^k.
\end{equation}
By Theorem~\ref{testar}, we identify $T^*\Pi E$ and $T^*\Pi E^*$. Rewritten in
coordinates $\x^i=(-1)^{\itt}\p^i,\p_i=\x_i$, the Hamiltonian $Q_{E^*}$ becomes
\begin{multline}\label{qestar}
  Q_{E^*}  =\p_i Q^{ia}(x) p_a+\frac{1}{2}\p_j\p_iQ^{ij}_k(x)(-1)^{\kt}\x^k=\\
  Q^{ia}(x) p_a\p_i+\frac{1}{2}(-1)^{\kt}Q^{ij}_k(x)\x^k\p_j\p_i.
\end{multline}
\begin{thm}[Roytenberg~\cite{roytenberg:thesis}, Vaintrob~\cite{vaintrob:algebroids}]\label{royt}
The Hamiltonians $S_E$ and $Q_{E^*}$ coincide.  Hence, $E$ is a Lie bialgebroid if and
only if the Hamiltonians $Q_E, Q_{E^*}\in \fun(T^*\Pi E)=\fun(T^*\Pi E^*)$ commute.
\end{thm}
\begin{proof}
The key statement is that $Q_{E^*}=S_E$ (after the identification
of $T^*\Pi E$ and $T^*\Pi E^*$). The second statement would follow
from Theorem~\ref{qsravno0}. We have to show that the Schouten
bracket on $\Pi E$ induced by the algebroid structure in $E^*$ and
specified by $S_E$ is the same as given by the quadratic
Hamiltonian~\eqref{qestar}. To this end, we  check the
non-vanishing Schouten brackets for coordinates. We get
$\{\x^i,x^a\}_S=\{\t(e^i),x^a\}_S=(-1)^{\itt}a(e^i)x^a=(-1)^{\itt}Q^{ia}$
and
$\{\x^i,\x^j\}_S=\t([e^i,e^j]_{E^*})=\t((-1)^{\jt}Q^{ij}_ke^k)=(-1)^{\itt+\kt}Q^{ij}_k\x^k$,
by the definition of the Schouten bracket induced by an algebroid
structure and formulas~\eqref{brackinestar}. (Here $\t{\co}
\fun(M,E^*)\to \fun(\Pi E)$ is $\a=e^i\a_i(x)\mapsto
\t(\a)=\x^i\a_i(x)$.) This exactly coincides with the Schouten
brackets given by the Hamiltonian~\eqref{qestar}. Hence,
$S_E=Q_{E^*}$.
\end{proof}

It follows  that the condition of a Lie bialgebroid is self-dual.
Under the identification $T^*\Pi E=T^*\Pi E^*$ the $QS$-structure
on $\Pi E$ corresponds to the $QS$-structure on $\Pi E^*$ so that
the homological field and the bracket effectively exchange places.
Note, however, that this self-duality is an extra statement not
necessary for writing the condition of a Lie bialgebroid in the
Poisson form $\{Q,S\}=0$.  It does not make sense in a nonlinear
case   considered in Section~\ref{generalized}.

How to get an analog of Drinfeld's double for a bialgebroid $E$? Roytenberg suggested the
following argument.  Consider the Hamiltonian vector fields on the supermanifold $T^*\Pi
E$ corresponding to the Hamiltonians $Q=Q_E$ and $S=Q_{E^*}$. According to the general
formulas, for an odd function on $T^*\Pi E$ we have
\begin{equation}
  X_f=\{f,\ \}=\der{f}{p_a}\der{}{x^a}-
  (-1)^{\at}\der{f}{x^a}\der{}{p_a}+\der{f}{\p_i}\der{}{\x^i}
  +(-1)^{\itt}\der{f}{\x^i}\der{}{\p_i}.
\end{equation}
Hence, we get:
\begin{multline}
  X_Q  = \x^iQ_i^a\der{}{x^a}-(-1)^{\at}\left(\der{Q_i^b}{x^a}p_b\x_i+
  \frac{1}{2}\der{Q_{ij}^k}{x^a}\p_k\x^j\x^i\right)\der{}{p_a} \\
  +\frac{1}{2}\x^j\x^iQ_{ij}^k\der{}{\x^k}
  +(-1)^{\kt}\bigl(Q_k^ap_a+\x^jQ_{jk}^i\p_i\bigr)\der{}{\p_k},
\end{multline}
\begin{multline}
  X_S  =\p_iQ^{ia}\der{}{x^a}-(-1)^{\at}\left(\der{Q^{ib}}{x^a}p_b\p_i
  +\frac{1}{2}(-1)^{\kt}\der{Q^{ij}_k}{x^a}\x^k\p_j\p_i\right)\der{}{p_a} \\
  +\bigl(Q^{ka}p_a+(-1)^{\itt}\p_jQ^{jk}_i\x^i\bigr)\der{}{\x^k}
  +(-1)^{\kt}\frac{1}{2}\p_j\p_iQ^{ij}_k\der{}{\p_k}.
\end{multline}
Mimicking the construction of Drinfeld's double for Lie bialgebras, define $\hat
Q_D:=X_Q+X_S$, and for the vector field $\hat Q_D$ on $T^*\Pi E$ we get
\begin{multline}
  \hat Q_D  = \bigl(\x^iQ_i^a+\p_iQ^{ia}\bigr)\,\der{}{x^a}\\
  -(-1)^{\at}\left(\der{Q_i^b}{x^a}p_b\x_i+\der{Q^{ib}}{x^a}p_b\p_i+
  \frac{1}{2}\der{Q_{ij}^k}{x^a}\p_k\x^j\x^i+
  \frac{1}{2}(-1)^{\kt}\der{Q^{ij}_k}{x^a}\x^k\p_j\p_i\right)\der{}{p_a} \\
  +\left(Q^{ka}p_a+(-1)^{\itt}\p_jQ^{jk}_i\x^i+\frac{1}{2}\x^j\x^iQ_{ij}^k\right)\der{}{\x^k}\\
  +(-1)^{\kt}\left(Q_k^ap_a+\x^jQ_{jk}^i\p_i+\frac{1}{2}\p_j\p_iQ^{ij}_k\right)\der{}{\p_k}.
\end{multline}
Obviously, the field $\hat Q_D$ is homological, because $Q$ and
$S$ commute. Some terms in this lengthy expression are
recognizable as similar to those in Drinfeld's double in the
bialgebra case (compare with~\eqref{qd}). However, because of the
presence of terms cubic in variables $\x^i,\p_j$, the homological
field $\hat Q_D$ cannot be related with a Lie algebroid over $M$.

In Roytenberg's thesis~\cite{roytenberg:thesis} it was suggested
to call the field $\hat Q_D\in\Vect(T^*\Pi E)$ the
\textit{Drinfeld double of the Lie bialgebroid $E$}.

It should be once again emphasized that the supermanifold $T^*\Pi
E$ does not have the form $\Pi A$ for any vector bundle $A\to M$
and that the field $\hat Q_D$ contains terms not appropriate for
Lie (anti)algebroids over $M$.

Two other constructions were suggested as analogs of Drinfeld's
double. Liu, Weinstein and Xu~\cite{weinstein:liuxu} considered
the direct sum $E\oplus E^*$ as a vector bundle over $M$ and
endowed it with a structure of a \textit{Courant algebroid}, a
notion that they defined. Roughly, it is a vector bundle with a
bracket of sections which satisfies the Jacobi identity up to
certain anomaly expressible in terms of a non-degenerate inner
product (see~\cite{weinstein:liuxu},
also~\cite{roytenberg:thesis}). There is a  natural projection
$T^*\Pi E\to \Pi(E\oplus E^*)$.
Roytenberg~\cite{roytenberg:thesis} showed that the Courant
algebroid structure in $E\oplus E^*$ is obtained by exactly the
same procedure (described  in Section~\ref{bistructures} above)
that reproduces the Lie bracket on the double of a Lie bialgebra
in case the base manifold is a point.
K.~Mackenzie~\cite{mackenzie:drinfeld,mackenzie:notions} suggested
to consider as a double of a Lie bialgebroid $E$ the manifold
$T^*E$ together with the structure of a ``double Lie algebroid'',
in a categorical sense precisely defined by him. The easy part of
this notion as is follows. There is a diagram
\begin{equation}\label{kirillsdouble}
\begin{CD}
    T^*E=T^*E^*@>>>E^*\\
    @VVV        @VVV\\
    E@>>>M
\end{CD}
\end{equation}
It can be seen that each arrow is a vector bundle projection, and
the horizontal arrows constitute a morphism of the vertical bundle
structure, and vice versa. Such structure is called a
\textit{double vector bundle}. Moreover, it turns out that there
are Lie algebroid structures on all four sides, and, in addition,
certain conditions are
satisfied~\cite{mackenzie:drinfeld,mackenzie:notions} which allow
to describe the diagram~\eqref{kirillsdouble} as a \textit{double
Lie algebroid} over $M$ (this is the hard part). Again, the
``direct sum object'' $E\oplus E^*$ sits inside this picture, due
to the projection $T^*E\to E\oplus E^*$. (I hope to elucidate the
relation of Mackenzie's beautiful picture with the super approach
elsewhere.)

Since there is no hope to interpret the homological field $\hat
Q_D$ in terms of a vector bundle over $M$, other approach should
be taken. We think that the key hint is in the following fact. One
can check that a change of coordinates in $T^*\Pi E$ (see the
Appendix) preserves the following $\ZZ$-gradings: $\hash p+\hash
\p$ and $\hash p+ \hash \x$ (the symbol $\hash$ means total degree
in a given variable), as well as their linear combinations, viz.,
$\hash \x-\hash\p$ (cf. with the similar gradings considered
in~\cite{tv:ivb,tv:class}.) It was observed
in~\cite{roytenberg:thesis} that the vector field $\hat Q_D$ is
homogeneous w.r.t. a grading defined as
\begin{equation}\label{weight}
  (\hash p+\hash \p)+(\hash p+ \hash \x)=2\hash p+\hash \p+\hash \x,
\end{equation}
of degree  $+1$.  For a vector bundle $E\to M$, one can uniquely
characterize the homological fields on $\Pi E$ of the
form~\eqref{algebroid} corresponding to a Lie algebroid structures
in $E$ as having the degree $+1$ in linear fiber coordinates
$\x^i$.

Hence,  \textit{a $\ZZ$-grading should be regarded as a
replacement for an absent linear structure}. The trick is that,
differently from the vector bundle case, one should allow to count
degrees in different variables with possibly different weights (a
situation well known in algebraic topology). We systematically
develop this approach in the next sections.

\section{Graded manifolds}\label{graded}

\begin{de}\label{defgraded}
Consider a domain $U\subset \mathbb R^{n|m}$ with coordinates
$x^a$ (some even, some odd). Assign to them \textit{weights}:
$w(x^a)=w_a\in \mathbb Z$. Consider the algebra of all smooth
functions of $x^a$ that are polynomial in coordinates of nonzero
weight. Denote it $A$. It is a $\mathbb Z$-graded algebra (as well
as $\Z$-graded; these two gradings are independent). We assume
that $U$ is cylindrical in the directions of even variables of
nonzero weight. There are naturally defined morphisms of pairs
$(U,A)$ (homomorphisms of algebras in particular  must respect
weight). A \textit{graded manifold} is a locally ringed space
modeled on such pairs. \textit{Morphisms} of graded manifolds
preserve weight of functions.
\end{de}

Simply speaking, a graded manifold $M$ is a supermanifold with a privileged class of atlases in
which particular coordinates are assigned numbers called ``weight'' so that the changes of
coordinates respect these weights.

On every graded manifold there is an action of the \textit{scaling transformations} (group $\RR^*$
under multiplication). In coordinates, $(x^a)\mapsto (\lambda^{w_a}\, x^a)$, $\lambda\in \RR^*$.

\begin{rem}
It is possible to consider non-integer weights as well. The difference between coordinates of zero
and nonzero weight is similar to the distinction between physical quantities that in some system of
units are dimensionless and those that are not (e.g., measured in centimeters). It is possible to
consider a function like sine or logarithm of a dimensionless quantity, but not of the one that is
measured in centimeters or grams, for which only homogeneous algebraic operations are permitted.
\end{rem}

{\footnotesize The idea of a ``graded manifold'' is, of course,
very natural. Formal ``$\ZZ$-graded manifolds'' were used by
Kontsevich in~\cite{kontsevich:quant}, however, his $\ZZ$-grading
simply underlies  $\Z$-grading. Our definition is very close to
the notion of ``$N$-manifolds'' suggested by \v{S}evera
in~\cite{severa:luminy}. However, \v{S}evera also couples
$\Z$-grading and $\ZZ$-grading.  The  point in
Definition~\ref{defgraded} is that parity and weight should be
considered completely independently.

}

\smallskip
\begin{exs} 1. An arbitrary supermanifold: assign zero weight to all coordinates.
2. A vector space: assign the same (arbitrarily chosen) number $N$
as weight to all linear coordinates. 3. A vector bundle: assign
$N$ for all linear fiber coordinates (even or odd) and assign zero
to all base coordinates. 4. A double vector bundle, like $TE$ or
$T^*E$ for a vector bundle $E$.
\end{exs}

For vector spaces and vector bundles  $N=1$ is a standard choice.
In general, for a vector bundle $E$ we can write $E(N)$ for the
choice $N$, so $E=E(1)$. Denote $E(N)^*:=E^*(-N)$. Then the
natural pairing has zero weight.

\begin{thm} Suppose $M$ is a graded manifold. Then the manifolds $TM$, $T^*M$, $TM\Pi$,
$\Pi T^*M$ are bigraded. One grading is the standard vector bundle grading, the other is the
induced weight.
\end{thm}
\begin{proof}
We can lift the scaling transformations to the bundles $TM$, $T^*M$, $TM\Pi$, $\Pi T^*M$. If
on $M$ we have $w(x^a)=w_a$, then for the natural coordinates in the fibers of $TM$, $T^*M$,
$TM\Pi$, $\Pi T^*M$ we obtain, respectively, $w(v^a)=w_a$, $w(p_a)=-w_a$, $w(dx^a)=w_a$,
$w(x^*_a)=-w_a$. It can be checked directly that coordinate changes respect these weights.
\end{proof}

Various total gradings can be constructed as linear combinations of these two.

Notice that on graded manifolds all geometric objects assume
weights. In particular, for a vector field $X\in \Vect(M)$, if
$X=X^a(x)\partial_a$, then $w(X)=w(X^a)-w_a$. Thus the induced
weight on $T^*M$ of the corresponding Hamiltonian $p(X)=X^ap_a$
equals the weight of $X$ on $M$. The canonical Poisson bracket on
$T^*M$ has induced weight zero, because the derivatives w.r.t.
$x^a$ and $p_a$ have opposite weights, which mutually cancel. That
means that $w(\{f,g\})=w(f)+w(g)$. It follows that the weight of
the Hamiltonian vector field $X_f$ on $T^*M$ is the same as that
of a function $f$. On the other hand, the fiberwise degree of the
bracket on $T^*M$ is $-1$. Similar facts hold for $\Pi T^*M$ and
the canonical Schouten bracket.

\begin{rem}
The topological structure of graded manifolds is an interesting question. Precisely: what kind
of restrictions for a (super)manifold poses the existence of atlases with coordinate changes
preserving weight.

First, it is rather obvious that if all weights are nonnegative
and no ``external constants'' of nonzero weight (see below) are
allowed, then under mild assumptions, every graded  manifold can
be decomposed into a tower of fibrations of the form
\begin{equation}
  M=M_N\to M_{N-1}\dots\to M_2\to M_1 \to M_0
\end{equation}
where coordinates on $M_0$ have zero weight, $M_1 \to M_0$ is a vector bundle, and all other
fibrations $M_{k+1}\to M_k$ are ``affine bundles'', in the sense that changes of variables are
linear in fiber coordinates plus additive terms of appropriate weight. (We group all
coordinates by decreasing weights $w_N>w_{N-1}>\ldots>w_2>w_1>w_0=0$;   at step $N$
coordinates in the fiber have weight $w_N$ and coordinates on the base $M_{N-1}$ have smaller
weights, etc.)

Second, similarly to supermanifolds (for which we need external
odd parameters or ``odd constants'' to be able to consider
universal families, etc.), for graded manifolds we might need, in
addition, to use external parameters (``constants'') with assigned
nonzero weights. They can enter formulas that would otherwise be
not homogeneous. The role of such constants is similar to
``universal constants'' in physics or simply to a choice of
physical units. An example of situation where this is necessary is
constructing metric (of zero weight) on the total space of a
vector bundle, from a metric on the base and a fiberwise metric in
the bundle (which has  weight $-2$ for the standard choice).
\end{rem}

\section{Generalized Lie bialgebroids and their doubles}\label{generalized}

Consider a graded $QS$-manifold $M$. Suppose that the homological field $\hat Q$ has weight $q$ and
the Schouten bracket $\{\ ,\ \}_S$ has weight $s$. To them correspond odd Hamiltonians (functions
on $T^*M$) $Q=p_{\hat Q}$ (fiberwise degree $1$) and $S$ (fiberwise degree $2$) of weights $q$ and
$s$, respectively. To functions $Q$ and $S$ correspond Hamiltonian vector fields $X_{Q}$ and $X_S$
on $T^*M$. It is convenient to arrange the information about gradings  in a table:

\par\medskip\noindent
\begin{center}
{\renewcommand{\arraystretch}{1.3}
\begin{tabular}{|c|c|c|c|c|c|c|c|} \hline
   \rule{0pt}{18pt}  & $\hat Q$ & $\{\ ,\ \}_S$& $Q$& $S$ & $\{\ ,\ \}$ & $X_{Q}$ & $X_S$\\ \hline
    weight on $M$ & $q$ & $s$ & -- & -- & -- & -- & --\\ \hline
    $w$ = induced weight on $T^*M$& -- & -- & $q$ &$s$& $0$&$q$&$s$\\ \hline
    $d$ = fiberwise degree& -- & -- & $1$ &$2$& $-1$&$0$&$1$\\ \hline
    $w+(q-s)d$& -- & -- & $2q-s$  & $2q-s$& $-q+s$ &$q$&$q$\\ \hline
\end{tabular}
}
\end{center}

\medskip\noindent
The last row follows by a direct calculation.

\begin{thm}\label{mydouble}
Define the total weight on $T^*M$ as $w+(q-s)d$, where $w$ is the
induced weight and $d$ the fiberwise degree. Define the vector
field $\hat Q_D=X_{Q}+\lambda X_S\in\Vect(T^*M)$ (where
$\lambda\in\RR$ is arbitrary, so it is a pencil of fields). Then
$T^*M$ with $\hat Q_D$ is a graded $Q$-manifold;  the homological
vector field $\hat Q_D$  is of weight $q$,  the same as the weight
of $\hat Q$ on $M$.
\end{thm}
\begin{proof}
Consider the Hamiltonian $Q_D=Q+\lambda S$. (Notice that it is odd, as it should be.) We have: $\{Q
,Q \}=p([\hat Q,\hat Q])=0$ (as $\hat Q^2=0$), $\{S ,S\}=0$ (because $S$ defines a Schouten
bracket), and $\{Q,S\}=L_{\hat Q} S=0$ (because $\hat Q$ is a derivation of the bracket). Thus
$\{Q_D,Q_D\}=0$. Hence, $[\hat Q_D,\hat Q_D]=2\hat Q_D^2=0$. The part concerning weights follows
from the table above.
\end{proof}

Suppose $\hat Q=Q^a(x)\lder{}{x^a}$ and $S=(1/2)S^{ab}(x)p_bp_a$. Then the explicit formula for
$\hat Q_D$ is
\begin{equation}\label{qdouble}
  \hat Q_D=\left(Q^a+p_bS^{ba}\vphantom{\der{S^{bc}}{x^a}}\right)\der{}{x^a}-
  (-1)^{\at}\left(\der{Q^b}{x^a}\,p_b+\frac{1}{2}\der{S^{bc}}{x^a}\,p_cp_b\right)\der{}{p_a}
\end{equation}
(we set for simplicity $\lambda=1$).

We denote $T^*M$ with the described graded $Q$-manifold structure
by $DM$ and call it the \textit{double} of a graded $QS$-manifold
$M$.  We can consider a graded $QS$-manifold $M$ as a generalized
Lie bialgebroid (more precisely, it is a generalized
antibialgebroid, i.e., a generalization of $\Pi E$ for a
bialgebroid $E$). The case $q=w(\hat Q)=+1$, $s=w(S)=-1$ is the
closest to the usual Lie bialgebroids.

\begin{ex}
Let $M=\Pi E$ for a Lie bialgebroid $E\to X$. Obviously, $DM$ as a
$Q$-manifold coincides with Roytenberg's ``Drinfeld's double'' of
the bialgebroid $E$. Let us calculate weights according to the
recipe of Theorem~\ref{mydouble}. For $M=\Pi E$ we have the
standard vector bundle weights: $w(x^a):=0$, $w(\x^i):=1$ (in the
notation of Section~\ref{bialgebroids}). The homological field has
weight $1$, the Schouten bracket weight $-1$
(see~(\ref{bialgebroid1}-\ref{bialgebroid2})). Then on $DM=T^*\Pi
E$ we get the induced weights $w(p_a)=0$, $w(\p_i)=-1$ and  the
fiberwise degree $d(p_a)=1$, $d(\p_i)=1$. Hence, our general
recipe for the total weight yields
$w+(1-(-1))d=w+2d=\#\x-\#\p+2(\#p+\#\p)=\#\x+\#\p+2\#p$, which is
exactly~\eqref{weight}.
\end{ex}

The construction of the double \textit{mutatis mutandis} carries over to  $QP$-manifolds. Let $M$
be a graded $QP$-manifold with the homological field $\hat Q$ of  weight $q$ and the Poisson
bracket $\{\ ,\ \}_P$ of weight $p$. To them correspond even multivector fields $Q'=\t({\hat
Q})=-Q^a(x)x^*_a$ (fiberwise degree $1$) and $P=(1/2)P^{ab}(x)x^*_bx^*_a$ (fiberwise degree $2$)
of weights $q$ and $p$, respectively. Consider the supermanifold $\Pi T^*M$ endowed with the
canonical Schouten bracket. We can again draw a table:

\par\medskip\noindent
\begin{center}
{\renewcommand{\arraystretch}{1.3}
\begin{tabular}{|c|c|c|c|c|c|c|c|} \hline
   \rule{0pt}{18pt}  & $\hat Q$ & $\{\ ,\ \}_P$& $Q'$& $P$ & $\{\ ,\ \}$ & $X_{Q'}$ & $X_P$\\ \hline
    weight on $M$ & $q$ & $p$ & -- & -- & -- & -- & --\\ \hline
    $w$ = induced weight on $\Pi T^*M$& -- & -- & $q$ &$p$& $0$&$q$&$p$\\ \hline
    $d$ = fiberwise degree& -- & -- & $1$ &$2$& $-1$&$0$&$1$\\ \hline
    $w+(q-p)d$& -- & -- & $2q-p$  & $2q-p$& $-q+p$ &$q$&$q$\\ \hline
\end{tabular}
}
\end{center}

\medskip\noindent
Here $X_{Q'}:=\{Q',\ \}$ and $X_P:=\{P,\ \}$ stand for the odd Hamiltonian vector fields
corresponding to the even functions $Q',P\in\fun(\Pi T^*M)$.

\begin{thm}\label{myodddouble}
Define the total weight on $\Pi T^*M$ as $w+(q-p)d$, where $w$ is
the induced weight and $d$ the fiberwise degree. Define the vector
field $\hat Q_D=X_{Q'}+\lambda X_P\in\Vect(\Pi T^*M)$ (where
$\lambda\in\RR$ is arbitrary). Then $\Pi T^*M$ with $\hat Q_D$ is
a graded $Q$-manifold, and  the homological vector field $\hat
Q_D$ is of weight $q$, which is the same as the weight of  $\hat
Q$ on $M$.
\end{thm}
\begin{proof}
Similar to the proof of Theorem~\ref{mydouble}. Consider the even function $Q_D'=Q'+\lambda P\in
\fun(\Pi T^*M)$. We have: $\{Q' ,Q' \}=\t([\hat Q,\hat Q])=0$ (as $\hat Q^2=0$), $\{P ,P\}=0$
(because $P$ defines a Poisson bracket), and $\{Q',P\}=L_{\hat Q} P=0$ (because it is equivalent to
$\hat Q$ being a derivation of the bracket $\{\ ,\ \}_P$, by the virtue of the non-degeneracy of
the canonical Schouten bracket, cf. the proof of Theorems~\ref{qsstructure} and~\ref{royt}). Thus
$\{Q_D',Q_D'\}=0$. Hence, $[\hat Q_D,\hat Q_D]=2\hat Q_D^2=0$. The part concerning weights follows
from the table.
\end{proof}

We call the graded $Q$-manifold $\Pi T^*M$ the \textit{double} of the graded $QP$-manifold $M$ and
denote it $DM:=\Pi T^*M$.

\smallskip
The double of a graded $QS$-manifold is automatically a
$QP$-manifold (with respect to the canonical symplectic
structure). In the same way,  the double of a graded $QP$-manifold
is automatically a $QS$-manifold, with respect to the canonical
odd symplectic structure. This corresponds to the existence of a
natural $\ad$-invariant inner product on Drinfeld's double
$\mathfrak d(\mathfrak g)=\mathfrak g\oplus \mathfrak g^*$. The
following theorem gives partial analogs of Drinfeld's and Manin's
theorems.

\begin{thm}
Let $M$ be a graded $QS$-manifold with the homological field $\hat Q$ and the Schouten bracket $\{\
,\ \}_S$. The homological vector field $\hat Q_D\in \Vect(DM)$ on the double $DM=T^*M$ has the
following properties:
\begin{enumerate}
  \item The canonical symplectic structure on $T^*M$ is $\hat Q_D$-invariant,
  \item $\hat Q_D$ is tangent to $M\subset DM$, and $\hat Q_{D{\displaystyle |M}}=\hat Q$,
  \item For arbitrary functions $f,g\in\fun(M)$ the Poisson bracket $\{f,\hat Q_D g\}$ also belongs to
  $\fun(M)$, and $\{f,\hat Q_D g\}=\{f,g\}_S$,
\end{enumerate}
and it is uniquely defined by these properties. Conversely, if for a graded manifold $M$, on
the cotangent bundle $T^*M$ there is a homological vector field $\hat Z$ such that it
preserves the canonical Poisson bracket, is tangent to $M\subset T^*M$, and for arbitrary
functions $f,g\in\fun(M)$ the bracket $\{f,\hat Z g\}$ is also a function on $M$, then $M$ is
a graded $QS$-manifold w.r.t. the vector field $\hat Q:=\hat Z_{{\displaystyle |M}}$ and the
odd bracket defined as $\{f,g\}_S:=\{f,\hat Z g\}$.
\end{thm}
\begin{proof}
Consider $\hat Q_D$ on the double $DM$. The statement (1) is clear. The statement (2) immediately
follows from the coordinate formula~\eqref{qdouble}. To obtain (3), notice that $\hat Q_Dg=\hat
Qg+\{S,g\}$, hence, because $\hat Qg\in\fun(M)$, it follows that $\{f,\hat
Q_D\}=\{f,\{S,g\}\}=\{f,g\}_S$, for arbitrary functions $f,g\in\fun(M)$). Now suppose that for an
arbitrary homological field $\hat Z\in\Vect(DM)$ the properties (1),(2),(3) hold. Then, by (1), we
conclude that $\hat Z$ is Hamiltonian, with some odd Hamiltonian $Z\in\fun(DM)$, where
\begin{equation} \label{Z}
  Z=Z_0(x)+Z^a(x)p_a+\frac{1}{2}\,Z^{ab}(x)\,p_bp_a+\frac{1}{3!}\,Z^{abc}(x)\,p_cp_bp_a\ldots,
\end{equation}
hence
\begin{multline}
  \hat Z=
  \left(Z^a(x)+p_bZ^{ba}(x)+\frac{1}{2}\,p_bp_cZ^{cba}(x)+\ldots\right)\der{}{x^a}
  \\-
  (-1)^{\at}\left(\der{Z_0}{x^a}+\der{Z^b}{x^a}\,p_b+
  \frac{1}{2}\,\der{Z^{bc}}{x^a}p_cp_b+\ldots\right)\der{}{p_a}.
\end{multline}
From the condition~(2) we conclude that $\lder{Z_0}{x^a}=0$, hence
w.l.o.g. $Z_0$ can be set to zero, and that $Z^a=Q^a$. Taking
$\{f,\hat Zg\}=\{f,\{Z,g\}\}$ for arbitrary functions on $M$ and
applying the condition~(3), we deduce that no terms of order $>2$
can be present in~\eqref{Z} and that $Z^{ba}=S^{ba}$. This proves
the uniqueness of $\hat Q_D$. Similar argument is applicable to an
arbitrary homological field $\hat Z\in\fun(T^*M)$ with the stated
properties  (without assuming a given $QS$-structure on $M$). We
get a field $\hat Q:=\hat Z_{{\displaystyle |M}}\in\Vect(M)$ and a
bracket $\{f,g\}_S:=\{f,\hat Z g\}$. That they make a
$QS$-structure follows directly.
\end{proof}

Similar theorem holds for  $QP$-manifolds.

\smallskip
Now, following the example of  Lie bialgebras, we want to provide the doubles of $QS$- and
$QP$-manifolds with brackets of the {\em same} parity as the original one. E.g., for a
$QS$-manifold $M$, we want a $\hat Q_D$-invariant Schouten bracket on $DM=T^*M$. To this end,
recall the construction of the Schouten bracket on $T^*\Pi\mathfrak g$  (see
Theorem~\ref{koskobkavduble}). We shall try to mimick it.

Let $x^a$ be local coordinates on $M$. Let us denote  the corresponding fiber coordinates in
$T^*M$ by $y_a$ now (the letter $p$ is reserved for a different use).  Consider the second
cotangent bundle $T^*T^*M$. Let $p_a$ and $q^a$ stand for fiber coordinates in $T^*(T^*M)$
(the conjugate momenta for $x^a$ and $y_a$, respectively). The changes of coordinates have the
following form (cf. equation~\eqref{cotanestar} in Appendix):
\begin{equation}
  \left\{\begin{aligned}
            x^a&=x^a(x')\\
            y_a&=J_{a}{}^{a'}(x')\,y_{a'}\\
            p_a&=\der{x^{a'}}{x^a}\,p_{a'}-(-1)^{\at\bt+\ct'}q^{c'}
            J_{c'}{}^{b}\der{J_{b}{}^{a'}}{x^a}\,y_{a'}\\
            (-1)^{\at}q^a&=(-1)^{\at'}q^{a'}J_{a'}{}^{a}
         \end{aligned}\right.
\end{equation}
Here $J_a{}^{a'}$ stands for the Jacobi matrix. Thus, $p_aq^a=(-1)^{\at}q^ap_a$ is not invariant.
However, if we are given a linear connection on $M$, so that for a local section of $T^*M$ the
covariant derivative is
\begin{equation}
  \nabla_a\a_b=\d_a\a_b-\G_{ab}{}^c\a_c,
\end{equation}
then the following ``long'' momentum
\begin{equation}
  P_a:=p_a+(-1)^{\bt(\at+1)}q^b\G_{ab}{}^cy_c=p_a+\G_{ab}{}^cy_cq^b
\end{equation}
transforms as
\begin{equation}
  P_a=\der{x^{a'}}{x^a}\,P_{a'},
\end{equation}
by a direct computation. Hence the expression
$P_aq^a\in\fun(T^*T^*M)$ is invariant. Denote it by $r$:
\begin{equation}
  r=P_aq^a=p_aq^a+\G_{ab}{}^cy_cq^bq^a.
\end{equation}
Following the example of Drinfeld's double for Lie bialgebras, we want to define the Schouten
tensor for $DM$ as $S_D:=(1/2)\{Q_D,r\}$.

In the following \textit{we assume that there is a connection of
zero weight} at our disposal. Then $w(\G_{ab}^c)=w_c-w_a-w_b$.
Hence, $w(\G_{ab}{}^cy_cq^b)=w_c-w_a-w_b-w_c+w_b=-w_a=w(p_a)$ ($w$
everywhere stands for weight on $M$ and the induced weight on
natural bundles), and the long momentum $P_a$ is homogeneous of
weight $-w_a$. We shall also use the  total weight  $W:=w+(q-s)d$
on $DM$ and the corresponding induced weight on $T^*DM$. Notice
that $W(x^a)=w_a$, $W(y_a)=-w_a+q-s$, $W(p_a)=-w_a$, $W(q^a)=
w_a-q+s$.

\begin{prop}
The function $S_D=(1/2)\{Q_D,r\}\in\fun(T^*DM)$ is quadratic in momenta $p_a, q^a$ and has weight
$W(S_D)=s$.
\end{prop}
\begin{proof}
Because $r$ is quadratic, $Q_D$ is linear and the canonical bracket has degree $-1$, the
function $S_D$ is quadratic. Now, let us calculate weight. The weight of $Q_D$ is the same as
that of $\hat Q_D$, $W(Q_D)=q$, and the weight of $r$ is $W(r)=-w_a+w_a-q+s=-q+s$. Hence,
because the canonical bracket has zero weight, $W(S_D)=W(Q_D)+W(r)=q-q+s=s$.
\end{proof}

Thus, the odd Hamiltonian $S_D$ defines on the double $DM$ an ``almost'' Schouten bracket of
weight $s$ automatically compatible with $\hat Q_D$. We shall give an explicit formula for it.
To this end, we may consider $x^a, y_a, P_a, q^a$ as convenient coordinates on $T^*DM$ (though
they are not Darboux). For a function $f\in\fun(DM)$ introduce a ``covariant partial
derivative''
\begin{equation}
  \nabla_af:=\der{f}{x^a}+\G_{ab}^cy_c\der{f}{y_b}.
\end{equation}
It transforms as a component of a covector on $M$. Notice that the
partial derivative $\lder{f}{y_a}$ transforms as a component of a
vector. We also adopt convention (usual in tensor calculus) that
that any operator nabla applied to a quantity with tensor indices
automatically incorporates extra Christoffel symbols.

\begin{lm}\label{skobkasrnam}
For an arbitrary function $f\in\fun(DM)$, the Poisson bracket of $p(X_f)$ with $r=P_aq^a$ is
given by the formula
\begin{multline}
\{p(X_f),r\}=\\
(-1)^{\ft(\at+\bt)}\left(\nabla_a\nabla_bf
-(-1)^{\bt\ct}R_{acb}^ky_k\der{f}{y_c}\right)q^bq^a
-(-1)^{(\ft+1)(\at+\bt)}\dder{f}{y_a}{y_b}\,P_bP_a.
\end{multline}
Here $R_{acb}^k$ is the Riemann tensor; the outer $\nabla_a$ contains an extra Christoffel
symbol for the tensor index of the inner $\nabla_b$.
\end{lm}

Applying this lemma to the functions $Q=Q^a(x)y_a$ and $S=(1/2)S^{ab}(x)y_by_a$, we get
\begin{thm}\label{tsdnam}
A $\hat Q_D$-invariant almost Schouten bracket on the double $DM$ of a graded $QS$-manifold
$M$ is given by the Hamiltonian $S_D=(1/2)\{Q_D,r\}$ which has the following explicit
expression
\begin{multline}\label{sdnam}
S_D=
\frac{1}{2}(-1)^{\at+\bt}\Biggl( \Bigl(\nabla_a\nabla_b Q^k -(-1)^{\bt\ct+\kt(\ct+1)}R_{acb}^k
Q^c\Bigr)y_k \\
+\Bigl(\frac{1}{2}\nabla_a\nabla_bS^{kl}-(-1)^{\bt\ct+(\kt+\lt)(\kt+\ct)}R_{acb}^l S^{kc}
\Bigr)y_ly_k\Biggr)q^bq^a \\
-\frac{1}{2}S^{ab}(x)P_bP_a.
\end{multline}
Here $R_{acb}^k$ is the Riemann tensor and nabla denotes  the covariant derivative of tensor
fields on the graded manifold $M$, w.r.t. the chosen linear connection.
\end{thm}

I skip the proofs of Lemma~\ref{skobkasrnam} and Theorem~\ref{tsdnam}. %
Notice that these statements are nonlinear analogs  of
Lemma~\ref{lfbar} and formula~\eqref{sdnabialgebre} for Lie
bialgebras. The difference with the Lie bialgebra case is due to
two independent reasons. One is a possible non-flatness of $M$,
hence the necessity to use a connection and the occurrence of
curvature as a result. The other difference is more substantial.
In a flat situation  one can replace covariant derivatives by
ordinary partial derivatives. In Drinfeld's case, when $Q^k$ is
quadratic in coordinates, the first term in~\eqref{sdnam} would
simply recover $Q^k$. For a general nonlinear $Q^k$, it would not
be so. This is a hint that we should not expect, generally,  for
the almost Schouten bracket given by $S_D$ to obey Jacobi.

\begin{ex}
Consider $M=\R{1|3}$ as a graded manifold, with an even coordinate
$x$ and odd coordinates $\x^i$ where we assign weights as $w(x)=2$
and $w(\x^i)=1$. Consider
\begin{equation}\label{eqnonlinalg}
    \hat Q=\left(x\x^1+\x^1\x^2\x^3\right)\der{}{x} +
    \x^1\x^3\der{}{\x^1}+\left(x+\x^1\x^2\right)\der{}{\x^2}.
\end{equation}
One can directly check that $\hat Q^2=0$. The homological
field~\eqref{eqnonlinalg} specifies on $\R{3|1}$ a structure of an
$L_{\infty}$-algebra, which is in fact homotopy equivalent to a
differential Lie superalgebra.  We want to consider $M$ as a
nonlinear analog of a Lie algebroid. Notice that $w(\hat Q)=+1$.
Take the identically zero Schouten bracket on $M$. Hence we have a
graded $QS$-manifold which is an analog of a Lie algebra
considered as a bialgebra (with the zero cobracket). The double of
$M$ is $DM=T^*M$ with the Hamiltonian homological vector field of
weight $+1$
\begin{multline*}
    \hat Q_D=\left(x\x^1+\x^1\x^2\x^3\right)\der{}{x} +
    \x^1\x^3\der{}{\x^1}+\left(x+\x^1\x^2\right)\der{}{\x^2}- \\
    \left(\x^1y+\h_2\right)\der{}{y}+
    \Bigl(\left(x+\x^2\x^3\right)y+\x^3\h_1+\x^2\h_2\Bigr)\der{}{\h_1}-
    \\
    \left(\x^1\x^3y+\x^1\h_2\right)\der{}{\h_2}+
    \left(\x^1\x^2y-\x^1\h_1\right)\der{}{\h_3}\,,
\end{multline*}
where $y$ and $\h_i$ are the conjugate momenta for $x$ and $\x^i$.
The ``total weights'' are: $W(x)=2$, $W(\x^i)=1$, $W(y)=0$,
$W(\h_i)=1$ (we assumed that $w(S)=-1$).  $\hat Q_D$   is just the
Lie derivative along $\hat Q$. (In the linear case, this
corresponds to a semidirect product Lie algebra structure in
$\mathfrak g\oplus \mathfrak g^*$ given by the coadjoint
representation.) We can construct a $\hat Q_D$-invariant almost
Schouten bracket on $DM$ as outlined above. Let the momenta
conjugate to $x$, $\x^i$, $y$, $\h_i$ on  $T^*DM=T^*T^*M$ be
denoted as $p$, $\pi_i$, $q$, $\k^i$, respectively. Take a flat
connection on $M$ so that the covariant derivatives w.r.t. $x$,
$\x^i$ are simply the partial derivatives. Then $r=pq+\pi_i\k^i$
and  either directly calculating or using Theorem~\ref{tsdnam} we
obtain $S_D=(1/2)\{Q_D,r\}$ as
\begin{equation*}
    S_D=-y\k^1q-\left(\x^3y+\h_2\right)\k^2\k^1+\left(\x^2y-\h_1\right)\k^3\k^1
    -\x^1y\k^3\k^2,
\end{equation*}
which gives the following invariant almost Schouten brackets:
\begin{align*}
    \{y,\h_1\}_{S_D}&=y\\
     \{\h_1,\h_2\}_{S_D}&=-\x^3y-\h_2\\
     \{\h_1,\h_3\}_{S_D}&=-\x^3y-\h_1\\
     \{\h_2,\h_3\}_{S_D}&=-\x^1y.
\end{align*}
Notice that they are homogeneous of weight $-1$. As expected, the
only non-zero brackets are between the ``dual'' coordinates $y,
\h_i$ (in the Lie algebra case we would recover the Lie--Poisson
bracket on $\mathfrak g^*$). However, they depend on variables
$\x^i$ as parameters, due to the nonlinearity of the original
field~\eqref{eqnonlinalg}. The Jacobi identity fails, as
$\{S_D,S_D\}_{T^*T^*M}=2\left((y\x^3+\h_2)\k^1\k^2\k^3-yq\k^1\k^3\right)\neq
0$, by a direct calculation. This can be attributed either to a
bad choice of a connection or to the fact that another
manifestation of our graded $Q$-manifold $M$ is a strong homotopy
Lie algebra
--- hence it would be plausible to look for ``homotopy'' Schouten
brackets on $DM$. In fact, both considerations can be related. As
mentioned above, the $L_{\infty}$-algebra defined by $\hat Q$ is
homotopy equivalent to an ordinary Lie superalgebra; i.e., the
cubic term in $\hat Q$ can be killed by a nonlinear
transformation.  There is no problem to construct a genuine
bracket on the dual of a Lie superalgebra; it should give a
``homotopy'' bracket on a homotopy equivalent object. Or, this
amounts to a different choice of a flat connection, in ``linear''
coordinates on $M$ corresponding to the Lie algebra structure
rather than in $x$, $\x^i$.
\end{ex}

The toy example above illustrates the general situation (except
that in general we do not have a ``linear'' reference object). It
points to a link between the non-linearity and homotopy . It is
probable that the double for a nonlinear $QS$- or $QP$-manifold
should be equipped with something like a ``homotopy Schouten
bracket'' rather than an ordinary one. This should be explored
further.

\section{Odd Lie bialgebras and the odd double}\label{odd}

In this section we return to the linear situation and specialize
our approach for odd Lie bialgebras, which were defined in
Section~\ref{bistructures}. Notice that an odd bialgebra structure
can be nontrivial (with a nontrivial odd cobracket) only in the
super case, for $\mathfrak g=\mathfrak g_0\oplus \mathfrak g_1$
where $\mathfrak g_1\neq 0$.

Let $\mathfrak g$ be an odd Lie bialgebra. We describe its structure using the
supermanifold $\Pi \mathfrak g$. Let the homological vector field $\hat
Q=(1/2)\x^j\x^iQ^k_{ij}\lder{}{\x^k}$  correspond to the Lie superalgebra structure in
$\mathfrak g$ and let the  Poisson tensor $P=(1/2)\x^k P^{ij}_k \x^*_j\x^*_i$ correspond
to the Lie superalgebra structure in $\Pi \mathfrak g^*$. Then the  Lie brackets in
$\mathfrak g$ are given by the formula
\begin{equation}\label{brack}
  [e_i,e_j]=(-1)^{\jt}Q^k_{ij}e_k,
\end{equation}
and the brackets in $\Pi \mathfrak g^*$ are given by
\begin{equation}\label{oddbrack}
  [\e^i,\e^j]=(-1)^{\itt+1}\e^k P^{ij}_k.
\end{equation}
Here $\e^i=\Pi e^i$ where the basis $e^i$  is right-dual to $e_j$.
The last formula follows from the definition of the Lie--Poisson
bracket (the Poisson bracket of the elements of $\Pi \mathfrak
g^*$ considered as linear functions on $\Pi \mathfrak g$ equals
their Lie bracket) and the formulas from Section~1.   It also
coincides with the one given by the duality argument, see
Example~\ref{linbivector} in the Appendix. (Notice that
$\x^i=\e^i$, as well as $x_i=e_i$.) It is not difficult to see
that the homological vector field on $\mathfrak g^*=\Pi\Pi
\mathfrak g^*$ defining the bracket~\eqref{oddbrack} is $\hat
Q^*=({1}/{2})(-1)^{\kt+1}P^{ij}_k x_jx_i\lder{}{x_k}$.

The odd bialgebra conditions for $\mathfrak g$ are
\begin{equation}\label{qp}
  \{Q,Q\}  =0, \quad
  \{P,P\} =0, \quad
  \{Q,P\}  =0,
\end{equation}
with the canonical Schouten brackets in $\Pi T^*\Pi \mathfrak g$, where $P$ is the
Poisson tensor and   $Q=\theta({\hat Q})=-(1/2)\x^j\x^iQ^k_{ij}\x^*_k\in \fun(\Pi T^*\Pi
\mathfrak g)$ is the even antivector field corresponding to the homological field $\hat
Q$. (There should be no confusion with  Hamiltonians used in the previous sections.)

\begin{lm} \label{nechkoskobka}
For every  Lie superalgebra $\mathfrak g$, on its antidual space
$\Pi\mathfrak g^*$ there is a naturally defined odd linear map
(``odd cobracket'')
\begin{equation}
  \delta{\co}  \Pi\mathfrak g^*\to S^2(\Pi\mathfrak g^*)\subset \Pi\mathfrak g^*\otimes
\Pi\mathfrak g^*.
\end{equation}
If the Lie bracket is given by~\eqref{brack}, then the odd
cobracket for it is given by the formula
\begin{equation}
  \delta(\e^k)=-\e^j\otimes \e^iQ_{ij}^k=-\e^j\e^iQ_{ij}^k.
\end{equation}
Here $\e_j\e_i=(1/2)(\e_j\otimes \e_i+(-1)^{(\itt+1)(\jt+1)}\e_i\otimes \e_j)$ is the
symmetric product.
\end{lm}
\begin{proof}
The map $\delta{\co}  \Pi\mathfrak g^*\to \Pi\mathfrak g^*\otimes \Pi\mathfrak g^*$ for
an arbitrary  Lie superalgebra $\mathfrak g$ is defined as follows. If we consider the
Lie bracket as the linear map $\delta^*{\co}  \mathfrak g\otimes \mathfrak g\to \mathfrak
g$, then the usual adjoint  gives the map $\delta'{\co}  \mathfrak g^* \to \mathfrak
g^*\otimes \mathfrak g^*$, and $\delta$ is the composition of $\delta'$ with the natural
isomorphism $\mathfrak g^*\otimes\mathfrak g^*\to \Pi \mathfrak g^*\otimes\Pi\mathfrak
g^*$ and the natural (odd) isomorphism $\Pi \mathfrak g^*\to\mathfrak g^*$ (left
multiplication by $\Pi$). In terms of bases we have:
\begin{equation*}
  \langle  {e_i\otimes e_j}, \delta'(e^k)\rangle
  =\langle [e_i, e_j], e^k\rangle=\langle (-1)^{\jt}Q_{ij}^l e_l,
  \e^k\rangle=(-1)^{\jt}Q_{ij}^k,
\end{equation*}
hence $\delta'(e^k)=(-1)^{\jt+\itt\jt}e^i\otimes e^j Q_{ij}^k=(-1)^{\itt+1}e^i\otimes e^j
Q_{ji}^k$. (Notice that for the basis $e^i$ right-dual to $e_i$,  $\langle {e_i\otimes
e_j},{e^k\otimes e^l}\rangle=(-1)^{\itt\jt}\delta_i{}^k\delta_j{}^l$ .) Now, we have
$\e^k\mapsto \Pi\e^k=\Pi\Pi e^k=e^k$, and we use the canonical identification
$\e^i\otimes \e^j= \Pi e^i\otimes \Pi e^j$ with $(-1)^{\itt}e^i\otimes  e^j$. Thus, for
the composition, we get $\delta{\co}  \e^k\mapsto e^k \mapsto (-1)^{\itt+1}e^i\otimes e^j
Q_{ji}^k\mapsto -\e^i\otimes\e^j Q_{ji}^k$. Or:
\begin{multline*}
  \delta(\e^k)=-\e^j\otimes\e^i Q_{ij}^k=-\frac{1}{2}
  \left(\e^j\otimes\e^i Q_{ij}^k+\e^i\otimes\e^j Q_{ji}^k\right)= \\
  -\frac{1}{2}
  \left(\e^j\otimes\e^i +\e^i\otimes\e^j (-1)^{(\itt+1)(\jt+1)}\right)Q_{ij}^k=
  -\e^j\e^i Q_{ij}^k.
\end{multline*}
In other words, $\delta=-2\hat Q$ restricted to linear functions
(with $\e^k=\x^k$).
\end{proof}

\begin{prop}
If $\mathfrak g$  is an odd bialgebra and the Lie bracket in $\Pi
\mathfrak g^*$ is given by formula~\eqref{oddbrack}, then  the odd
cobracket on $\mathfrak g$ is an odd linear map
\begin{equation}
  \delta{\co} \mathfrak g\to S^2(\mathfrak g)\subset\mathfrak g\otimes \mathfrak g
\end{equation}
defined by
the formula
\begin{equation}\label{oddcobrackg}
  \delta(e_k)=(-1)^{\kt}P^{ij}_ke_j\otimes e_i=(-1)^{\kt}P^{ij}_ke_je_i.
\end{equation}
\end{prop}
\begin{proof}
Immediately follows from~Lemma~\ref{nechkoskobka}. For an odd
bialgebra, to obtain formula~\eqref{oddcobrackg} for the odd
cobracket in $\mathfrak g$ we simply use the remark at the end of
proof for~Lemma~\ref{nechkoskobka} together with the above formula
for the vector field $\hat Q^*$. (Notice that $e_i=x_i$.)
\end{proof}

\begin{ex}\label{q1}
Consider a $1|1$-dimensional Lie superalgebra $\mathfrak g$ with a
basis $e\in\mathfrak g_0$, $\e\in\mathfrak g_1$ such that
$[\e,\e]=2e$ (the ``supersymmetry algebra''). There is an odd
inner product in $\mathfrak g$ defined by the condition
$(e,\e)=-(\e,e)=1$. One can check that this product is
$\ad\mathfrak g$-invariant. Moreover, we can use this product to
obtain a cobracket on $\mathfrak g$ (as odd-dual to the given Lie
bracket). We have
\begin{align}
  \delta(e) & =0, \label{koskobkanag1}\\
  \delta(\e) &= e\otimes e.\label{koskobkanag2}
\end{align}
The bracket structure in $\mathfrak g$ corresponds to the homological field $\hat
Q=-x^2\,\lder{}{\x}\in\Vect(\mathfrak g)$, and the odd cobracket corresponds to the
Poisson tensor $P=-x(\x^*)^2$. Clearly, the Schouten bracket $\{Q,P\}$ is zero. Hence,
the cobracket~\eqref{koskobkanag1},\eqref{koskobkanag2}  makes $\mathfrak g$ into an odd
bialgebra.
\end{ex}

Now let us elaborate the construction of the odd double briefly introduced in
Section~\ref{bistructures}. In the sequel we  use $\mathfrak d(\mathfrak g):=\mathfrak
g\oplus \Pi \mathfrak g^*$ as the notation for the odd double. The odd double has a
structure of an odd bialgebra defined as follows.

To get the Lie bracket in $\mathfrak d(\mathfrak g)$, consider the sum $Q+P\in \fun(\Pi
T^*\Pi \mathfrak g)$ of multivector fields $Q$ and $P$. By the equations~\eqref{qp}, it
satisfies $\{Q+P,Q+P\}=0$. Hence the corresponding odd Hamiltonian vector field $X_{Q+P}$
is homological. Calculating by the general formulas from Section~1, we get
\begin{align}
  X_Q & = -\frac{1}{2}\,\x^j\x^i Q_{ji}^k\der{}{\x^k}-(-1)^{\kt(\jt+1)}\x^*_j\x^i Q_{ik}^j\der{}{\x^*_k},\\
  X_P & = (-1)^{\jt}\x^*_j\x^i P^{jk}_i \der{}{\x^k} +
  (-1)^{\kt(\itt+\jt+1)}\frac{1}{2}\,\x^*_j\x^*_i P^{ij}_k \der{}{\x^*_k}.
\end{align}
Define the vector field $\hat Q_D$ on $\Pi T^*\Pi \mathfrak g$ as $\hat Q_D:=-X_{Q+P}$
(notice the minus sign). It is homological; explicitly:
\begin{multline}\label{oddalgebra}
  \hat Q_D=\frac{1}{2}\,\x^j\x^iQ^k_{ij}\der{}{\x^k}
            +\x^*_j\x^i\left(-(-1)^{\jt} P^{jk}_i \der{}{\x^k}+
            (-1)^{\kt(\jt+1)} Q_{ik}^j\der{}{\x^*_k}\right)-\\
            -(-1)^{\kt(\itt+\jt+1)}\frac{1}{2}\,\x^*_j\x^*_i P^{ij}_k \der{}{\x^*_k}.
\end{multline}
Conversely, the condition that $\hat Q_D^2=0$ is equivalent to the condition
$\{Q+P,Q+P\}=0$, i.e., to the odd Lie bialgebra conditions~\eqref{qp}. The field $\hat
Q_D$ corresponds to the Lie superalgebra structure in $\mathfrak d(\mathfrak g)$. Notice
that $\Pi T^*\Pi \mathfrak g\cong \Pi\mathfrak g\oplus \mathfrak g^*=\Pi(\mathfrak
g\oplus \Pi \mathfrak g^*)$. Here we can identify $x_i=\x^*_i$. If we translate the
vector field~\eqref{oddalgebra} into the language of Lie brackets, we get the following
formulas
\begin{align}
  { }[e_i,e_j]&=(-1)^{\jt}Q^k_{ij}e_k \label{brackinodddouble1} \\
  [e_i,\e^j]&=P_i^{jk}e_k + \e^kQ_{ki}^j \label{brackinodddouble2} \\
  [\e^i,\e^j]&=(-1)^{\itt+1}\e^k P^{ij}_k \label{brackinodddouble3}
\end{align}
in the odd double, by a straightforward calculation using
formula~\eqref{liebracket}). (Hence  $\mathfrak g\subset \mathfrak
d(\mathfrak g)$ and $\Pi\mathfrak g^*\subset\mathfrak d(\mathfrak
g)$ are subalgebras.)

To get the odd cobracket in $\mathfrak d(\mathfrak g)$,  we need a
Poisson bracket on the supermanifold $\Pi \mathfrak d(\mathfrak
g)=\Pi T^*\Pi \mathfrak g$. To this end, consider the
anticotangent bundle $\Pi T^*\Pi T^*\Pi \mathfrak g$. To avoid
confusion, on $\Pi T^*\Pi \mathfrak g$ rename $\x^*_i=:x_i$, so
the coordinates on $\Pi T^*\Pi \mathfrak g$ will be $\x^i, x_i$.
Denote now by $\x^*_i$, $x^{*i}$ the odd conjugate momenta for the
coordinates $\xi^i$, $x_i$ respectively. Clearly, $x^{*i}$
transforms as $\x^i$, hence the odd function
$\rho=x^{*i}\x^*_i=\x^*_ix^{*i}$  on $\Pi T^*\Pi T^*\Pi \mathfrak
g$ is well-defined. Consider the function $P_D$ defined as
$P_D=(1/2)L_{\hat Q_D}\rho=(1/2)\{-\t({\hat Q_D}),\rho\}$. It
follows that $P_D$ is a $\hat Q_D$-invariant even bivector field
on $\Pi T^*\Pi \mathfrak g$. To find it explicitly we can either
prove and apply an analog of Lemma~\ref{lfbar} or just calculate
straightforwardly.
\begin{prop}
\begin{equation}\label{pd}
  P_D=\frac{1}{2}\,x^{*j}x^{*i}Q_{ij}^k x_k+ \frac{1}{2}\,\x^kP_k^{ij}\x^*_j\x^*_i
\end{equation}
\end{prop}
\begin{proof}
Consider $\rho=\x^*_jx^{*j}$ and
\begin{multline}
  -\t(\hat Q_D)=\frac{1}{2}\,\x^j\x^iQ^k_{ij}\x^*_k
            +\x^*_j\x^i\left(-(-1)^{\jt} P^{jk}_i \x^*_k+
            (-1)^{\kt(\jt+1)} Q_{ik}^jx^{*k}\right)-\\
            -(-1)^{\kt(\itt+\jt+1)}\frac{1}{2}\,\x^*_j\x^*_i P^{ij}_k x^{*k}.
\end{multline}
Their Schouten bracket equals
\begin{multline*}
    \{-\t({\hat Q_D}),\rho\}=(-1)^{\jt}\,\der{(-\t({\hat
    Q_D}))}{\x^j}\der{\rho}{\x^*_j}-(-1)^{\jt}\,\der{(-\t({\hat
    Q_D}))}{x_j}\der{\rho}{x^{*j}}= \\
    (-1)^{\jt}\biggl(\x^iQ_{ij}^k\x^*_k + (-1)^{\itt(\jt+1)}x_i\left(-(-1)^{\itt}P_j^{ik}\x^*_k +
    (-1)^{\kt(\itt+1)}Q_{jk}^ix^{*k}\right)\biggr)x^{*j}\,- \\
    (-1)^{\jt}\biggl(\x^i\left(-(-1)^{\jt}P_i^{jk}\x^*_k + (-1)^{\kt(\jt+1)}Q_{ik}^jx^{*k}\right)
    -(-1)^{\kt(\itt+\jt+1)}x_iP_k^{ij}x^{*k}\biggr)\x^*_j= \\
    (-1)^{\jt}\x^iQ_{ij}^k\x^*_k x^{*j} - (-1)^{\jt+\kt(\jt+1)}Q_{ik}^jx^{*k}\x^*_j \\
    -(-1)^{\jt+\itt(\jt+1)+\itt}x_iP_j^{ik}\x^*_k x^{*j} +
    (-1)^{\jt+\kt(\itt+\jt+1)}x_iP_k^{ij}x^{*k}\x^*_j + \\
    (-1)^{\jt+\kt(\itt+1)+\itt(\jt+1)}x_i Q_{jk}^ix^{*k}x^{*j}
    + \x^iP_i^{jk}\x^*_k\x^*_j.
\end{multline*}
After exchanging the indices $j$ and $k$ in the second and the fourth terms of the last
expression, changing order of factors and calculating signs, it turns out that these
terms cancel with the first and the second terms respectively. In the same way, the
factors in the fifth term can be rearranged so that it completely cancels  signs.
Finally:
\begin{equation}
  \{-\t({\hat Q_D}),\rho\}=x^{*k}x^{*j}Q_{jk}^ix_i + \x^iP_i^{jk}\x^*_k\x^*_j,
\end{equation}
which differs from the required formula only by the notation for the indices.
\end{proof}
\begin{cor}
$P_D$ is a Poisson tensor on $\Pi \mathfrak d(\mathfrak g)$.
\end{cor}
\begin{proof} The two terms in~\eqref{pd} have the zero canonical Schouten bracket
because they do not contain conjugate
variables. Each term is a Poisson tensor by itself. Indeed, the
second term is the Poisson tensor on $\Pi\mathfrak g$. The first
term equals $P^*=(1/2)(-1)^{\itt+\jt}Q_{ij}^k x_kx^{*j}x^{*i}$,
the Poisson tensor that gives the Lie--Poisson bracket
$\{x_i,x_j\}=(-1)^{\jt}Q_{ij}^k x_k$ on $\mathfrak g^*$ generated
by the Lie bracket in $\mathfrak g$.
\end{proof}

\begin{cor}
The nonvanishing Poisson brackets of coordinates on $\Pi \mathfrak d(\mathfrak g)$
are:
\begin{equation}
  \{\x^i,\x^j\}_{P_D}=(-1)^{\itt+1}\x^k P^{ij}_k, \quad \{x_i,x_j\}_{P_D}=(-1)^{\jt}Q_{ij}^k x_k.
\end{equation}
\end{cor}

\begin{cor}
The odd cobracket $\dd{\co} \mathfrak d(\mathfrak g)\to S^2\,\mathfrak d(\mathfrak
g)\subset\mathfrak d(\mathfrak g)\otimes \mathfrak d(\mathfrak g)$ is given by the
formulas:
\begin{align}
  \delta(e_k)&=
                (-1)^{\kt}P^{ij}_ke_je_i, \label{cobrackinodddouble1}\\
  \delta(\e^k)&=
                -\e^j\e^iQ_{ij}^k. \label{cobrackinodddouble2}
\end{align}
Hence,  both $\mathfrak g\subset \mathfrak d(\mathfrak g)$ and $\Pi\mathfrak
g^*\subset\mathfrak d(\mathfrak g)$ are sub-bialgebras.
\end{cor}

It is possible  to formulate everything entirely in terms of linear spaces $\mathfrak g$
and $\Pi \mathfrak g^*$. In the corresponding analogs of Theorem~\ref{drinfeld} and
Theorem~\ref{manin}, there appear odd inner products.

\begin{thm}\label{odddrinf}
Let $\mathfrak g$ be an odd Lie bialgebra. Then in the vector space $\mathfrak
d(\mathfrak g):=\mathfrak g\oplus \Pi \mathfrak g^*$ there is a natural structure of an
odd Lie bialgebra uniquely characterized by the following properties:
\begin{enumerate}
  \item $\mathfrak g$ and $\Pi\mathfrak g^*$ are Lie subalgebras in $\mathfrak d(\mathfrak g)$,
  \item The natural odd inner product in $\mathfrak d(\mathfrak
g)=\mathfrak g\oplus \Pi \mathfrak g^*$ is $\ad$-invariant,
  \item The odd cobracket $\dd{\co} \mathfrak d(\mathfrak g)\to S^2\,
  \mathfrak d(\mathfrak g)\subset\mathfrak d(\mathfrak g)\otimes \mathfrak d(\mathfrak g)$
  is given by the coboundary of the (odd) invariant element
  $\rho=e_i \e^i\in S^2(\mathfrak g\oplus \Pi\mathfrak g^*)$.
  \end{enumerate}
\end{thm}
\begin{proof}
The explicit formulas defining the odd bialgebra structure in $\mathfrak d(\mathfrak g)$
are given above:~(\ref{brackinodddouble1}-\ref{brackinodddouble3})
and~(\ref{cobrackinodddouble1}-\ref{cobrackinodddouble2}). We shall check uniqueness and
the statement concerning the cobracket. Recall that an odd symmetric inner product has
the property
\begin{equation}
  (u,v)=(-1)^{\ut+\vt+\ut\vt}(v,u).
\end{equation}
In $\mathfrak g\oplus \Pi \mathfrak g^*$ the nonvanishing inner products of the basis
vectors are $(e_i,\e^j)=\dd_i{}^j,$ $(\e^j,e_i)=-\dd_i{}^j$. The $\ad$-invariance
condition for an odd inner product  means that
\begin{equation}
  ([u,v],w)+(-1)^{\ut(\vt+1)}(v,[u,w])=0
\end{equation}
for all $u,v,w$, or
\begin{equation}
  ([v,u],w)=(-1)^{\ut}(v,[u,w]).
\end{equation}
Suppose $[e_i,\e^j]=f_i^{jk}e_k+\e^k g_{ki}^j$ with indeterminate
coefficients $f_i^{jk}$, $g_{ki}^j$. Then
$([e_i,\e^j],\e^k)=f_i^{jk}$, $(e_k,[e_i,\e^j])=g_{ki}^j$. From
the $\ad$-invariance we get:
$([e_i,\e^j],\e^k)=(-1)^{\jt+1}(e_i,[\e^j,\e^k])=(-1)^{\jt+1}(e_i,(-1)^{\jt+1}\e^lP_l^{jk})=P_i^{jk}$,
and also
$(e_k,[e_i,\e^j])=(-1)^{\itt}([e_k,e_i],\e^j)=(-1)^{\itt}((-1)^{\itt}Q_{ki}^l,\e^j)=Q_{ki}^j$,
where we used formulas~\eqref{brackinodddouble1}
and~\eqref{brackinodddouble3} that
follow from the 
condition 1. Hence $f_i^{jk}=P_i^{jk}$, $g_{ki}^j=Q_{ki}^j$, and
we recover formula~\eqref{brackinodddouble2}. The uniqueness is
proved. Consider now the element
$\rho=e_i \e^i=(1/2)(e_i \otimes \e^i+\e^i \otimes e_i)$. 
Its coboundary in the Lie algebra cochain complex is the function
$d\rho(u)=(-1)^{\ut}(S^2\ad u) (\rho)$. Hence,
\begin{multline}
  d\rho(e_k)=(-1)^{\kt}\left([e_k,e_i]\e^i+(-1)^{\kt\itt}\e^i[e_k,e_i]\right)= \\
  (-1)^{\kt}\left((-1)^{\itt}Q_{ki}^l e_l\e^i+(-1)^{\kt\itt}\e^i(P_k^{il} e_l+\e^l Q_{lk}^i)\right)=
  (-1)^{\kt}P_k^{il} e_le_i,
\end{multline}
where we omitted a straightforward simplification. Similarly,
\begin{multline}
  d\rho(\e^k)=(-1)^{\kt+1}\left([\e^k,e_i]\e^i+(-1)^{(\kt+1)\itt}\e^i[\e^k,e_i]\right)=
  \\
  (-1)^{\kt+1}\left((-1)^{\kt}\e^i[e_i,\e^k]+(-1)^{(\kt+1)\itt}\e^i[\e^k,e_i]\right)= \\
  (-1)^{\kt+1}\left((-1)^{\kt}\e^i(P_i^{kl}e_l+\e^lQ_{li}^k)+
  (-1)^{(\kt+1)\itt}\e^i(-1)^{\kt+1}\e^lP_l^{ki}\right)=
  -\e^i\e^lQ_{li}^k.
\end{multline}
Thus, $d\rho=\dd$.
\end{proof}

\begin{thm}
Let the vector space $\mathfrak d=\mathfrak a \oplus \mathfrak b$
have a structure of a Lie superalgebra with an invariant odd inner
product. Suppose that the subspaces $\mathfrak a$ and $\mathfrak
b$ are isotropic subalgebras in $\mathfrak d$. Then $\mathfrak b
\cong \Pi\mathfrak  a^*$, $\mathfrak a \cong \Pi\mathfrak  b^*$,
and $\mathfrak a$, $\mathfrak b$ are odd Lie bialgebras which are
in odd duality. The Lie superalgebra structure on $\mathfrak d$ is
isomorphic to that of the odd double either of $\mathfrak a$ or of
$\mathfrak b$.
\end{thm}
\begin{proof}It is possible to pick bases $e_i\in\mathfrak a$ and $\e^i\in\mathfrak b$,
$\tilde e_i=\itt$, $\tilde \e^i=\tilde e_i+1=\itt+1$ such that
$(e_i,\e^j)=\dd_i{}^j$, $(\e^j,e_i)=-\dd_i{}^j$. By an argument as
in the proof of the previous theorem, we immediately recover the
commutators of the basis elements in the
form~(\ref{brackinodddouble1}--\ref{brackinodddouble3}). Because
we are given that $\mathfrak d$ is a Lie superalgebra, the vector
field $\hat Q_D$ defining the
brackets~(\ref{brackinodddouble1}--\ref{brackinodddouble3}) is
homological. Hence, its coefficients obey the relations that can
be written as $\{Q,Q\}=0$, $\{P,P\}=0$, $\{Q,P\}=0$, and we arrive
at an odd bialgebra structure in $\mathfrak a$. The rest is
obvious. Notice that $\mathfrak a$ and $\mathfrak b$ can be
interchanged.
\end{proof}

In the classical case of Drinfeld's  Lie bialgebras, one of the
immediate applications of Drinfeld's double and of the concept of
Manin's triples is a natural bialgebra structure in the classical
Lie algebras like $\mathfrak{gl}(n)$ or $\mathfrak{sl}(n)$ and,
more generally, in symmetrizable Kac--Moody Lie algebras. Roughly
speaking, these algebras turn out to be ``almost'' doubles, i.e,
are obtained from doubles by the factorization by a Lie bialgebra
ideal isomorphic to an ``extra copy'' of a Cartan subalgebra
(see~\cite{charipressley} for details). As it turns out, this does
not carry over straightforwardly to obvious candidates for an odd
Lie bialgebra such as the Lie superalgebra $\mathfrak q(n)$ (see
below). However, we can succeed by slightly generalizing  the
notion of the double, as follows.

Consider a Lie (super)algebra $\mathfrak a$ and  a Lie
(super)algebra $\mathfrak h$. We suppose the following. Let
$\mathfrak h$ act on $\mathfrak a$ by derivations, and let
$\mathfrak h$ be endowed with an invariant inner product of parity
$\a$. For $\a=0$ or $\a=1$ assume that there is a Lie bracket  in
$\mathfrak{a}^*$ or $\Pi\mathfrak{a}^*$, respectively. Let the
contragredient action of $\mathfrak h$ on $\mathfrak{a}^*$ or
$\Pi\mathfrak{a}^*$ be  also by derivations. Denote
$\mathfrak{b}:=\mathfrak{a}^*$ or
$\mathfrak{b}:=\Pi\mathfrak{a}^*$ as appropriate. Then
$\mathfrak{a}\oplus\mathfrak{h}\oplus\mathfrak{b}$ will be a
potential ``double    of $\mathfrak{a}$ (or
$\mathfrak{a}\oplus\mathfrak{h}$)  over $\mathfrak{h}$''. More
precisely, this vector space has a natural inner product of parity
$\a$. Require a Lie bracket on
$\mathfrak{a}\oplus\mathfrak{h}\oplus\mathfrak{b}$ such that both
$\mathfrak{a}\oplus\mathfrak{h}$ and
$\mathfrak{b}\oplus\mathfrak{h}$ (with natural brackets) are
subalgebras and the inner product is invariant. It is a
straightforward check that such conditions specify the bracket
uniquely (if it exists). This is completely similar, e.g.,  to
proofs of Drinfeld's theorem or Theorem~\ref{odddrinf}. The
resulting formula for the ``cross'' bracket would be
\begin{equation}
  [e_i,e^j]=\pm Q_i^{jk}e_k \pm e^kQ_{ki}^j
  \pm Q_{i\mu}^{j}g^{\mu\lambda}e_{\lambda};
\end{equation}
here we use the notation $e^j$ for the dual basis either in
$\mathfrak{a}^*$ or in $\Pi\mathfrak{a}^*$, and $[e_i,e_j]=\pm
Q_{ij}^ke_k$, $[e^i,e^j]=\pm e^k Q_k^{ij}$, $[e_{\mu},e_i]=\pm
Q_{i\mu}^k e_k$ (the last formula stands for the action of
$\mathfrak{h}$ on $\mathfrak{a}$), and $g^{\mu\l}$ is the inverse
Gram matrix for $\mathfrak{h}$. We call $\mathfrak{a}$ a
\textit{Lie bialgebra over $\mathfrak{h}$} (\textit{of parity
$\a$}) if in $\mathfrak{a}\oplus\mathfrak{h}\oplus\mathfrak{b}$ is
obtained a genuine Lie bracket. Consider an invariant element $e_i
e^i$ of parity $\a$, where the product means either the wedge
product or the symmetric product. Applying $\ad$ of elements of
$\mathfrak{a}\oplus\mathfrak{h}\oplus\mathfrak{b}$, we arrive at a
cobracket $\delta$ on
$\mathfrak{a}\oplus\mathfrak{h}\oplus\mathfrak{b}$ (even or odd).
One can see that $\delta$ vanishes on $\mathfrak{h}$ (because
$\mathfrak{h}$ acts, up to a sign, by the adjoint operators on
$e_i$ and $e^i$) and that
\begin{align}
  \delta(e_k) & =\pm Q_k^{ij}e_i e_j \pm Q_{k\mu}^{j}g^{\mu\l}e_j e_{\l}, \label{mycobr1}\\
  \delta(e^k) & =\pm Q^k_{ij}e^i e^j \pm Q^{k}_{j\mu} g^{\mu\l}e^j e_{\l}. \label{mycobr2}
\end{align}
Immediately follows that  cobrackets~\eqref{mycobr1} and
\eqref{mycobr2} separately satisfy co-Jacobi, because they are
nothing but the cobrackets dual to the Lie brackets in
$\mathfrak{b}\oplus\mathfrak{h}$ and
$\mathfrak{a}\oplus\mathfrak{h}$, respectively. Hence, together we
have a Lie bracket in the dual (or antidual, depending on $\a$)
space to
$\mathfrak{g}=\mathfrak{a}\oplus\mathfrak{h}\oplus\mathfrak{b}$.
By the construction, it is an $\ad \mathfrak{g}$-cocycle. Hence
$\mathfrak{g}=\mathfrak{a}\oplus\mathfrak{h}\oplus\mathfrak{b}$ is
a Lie bialgebra (even or odd, depending on $\a$), which is, by
definition, the \textit{double of $\mathfrak{a}\oplus
\mathfrak{h}$ over $\mathfrak h$}. This argument can be reversed,
leading to such a ``relative double'' structure in a Lie
(super)algebra with an invariant inner product and an appropriate
decomposition into
$\mathfrak{a}\oplus\mathfrak{h}\oplus\mathfrak{b}$. In particular,
the even case examples with Kac--Moody algebras
(see~\cite{charipressley}) are covered by this construction. For
us the main motivation is the following example (where we get an
odd bialgebra).

\begin{ex}
Consider the matrix Lie superalgebra $\mathfrak{q}(n)$. Recall that
\begin{equation}
    \mathfrak{q}(n)=\bigl\{x\in \mathfrak{gl}(n|n)\,\bigl|\, [x,I]=0 \bigr.\bigr\},
    \text{ where }
   I= \begin{pmatrix}
  0 & 1 \\
  1 & 0
\end{pmatrix}.
\end{equation}
It consists of even matrices of the appearance
\begin{equation}\label{evenq}
x=\begin{pmatrix}
  a & b \\
  -b & a
\end{pmatrix}
\end{equation}
and odd matrices of the appearance
\begin{equation}\label{oddq}
x=\begin{pmatrix}
  a & b \\
  b & -a
\end{pmatrix}.
\end{equation}
(For the numerical entries, $b=0$ in~\eqref{evenq} and $a=0$
in~\eqref{oddq}.) There is an invariant odd inner product on
$\mathfrak{q}(n)$, defined as
\begin{equation}
  (x,y)=(-1)^{\xt}\otr (xy),
\end{equation}
where the ``odd trace'' is an (essentially unique) trace on
$\mathfrak{q}(n)$ given by the formula
\begin{equation}
  \otr x=\frac{1}{2}\str (Ix)=\tr b.
\end{equation}
There is an analog of the Cartan decomposition:
$\mathfrak{q}(n)=\mathfrak{n}_{+}\oplus\mathfrak{h}\oplus\mathfrak{n}_{-}$
into the ``upper-triangular'', ``lower-triangular'' and
``diagonal'' matrices. The Cartan subalgebra $\mathfrak{h}\subset
\mathfrak{q}(n)$ is no longer commutative. It consists of $n$
copies of the $1|1$-dimensional superalgebra considered in
Example~\ref{q1}. (It is exactly the non-commutativity of
$\mathfrak{h}$ that prevents to define a bialgebra structure in
$\mathfrak{q}(n)$ similarly to $\mathfrak{gl}(n)$.) We can apply
the above considerations and interpret $\mathfrak{q}(n)$ as a
``relative'' double of $\mathfrak{n}_{+}\oplus\mathfrak{h}$ over
$\mathfrak{h}$. Thus $\mathfrak{q}(n)$ has a natural structure of
an odd bialgebra. Direct calculations of commutators yield the
following explicit formulas for the odd cobracket
\begin{equation}
  \dd{\co}  \mathfrak{q}(n)\to S^2\mathfrak{q}(n).
\end{equation}
For $k=l$:
\begin{align}
  \dd(e_{kk}) & =0, \quad \dd(\e_{kk})=0; \label{koskobkanakartane}\\
  \intertext{for $k<l$:}
  \dd(e_{kl}) &
  =2\sum_{k<i<l}\left(-e_{ki}\e_{il}+\e_{ki}e_{il}\right)+(\e_{kk}-\e_{ll})e_{kl}-
  (e_{kk}-e_{ll})\e_{kl},\\
  \dd(\e_{kl}) &
  =2\sum_{k<i<l}\left(e_{ki}e_{il}-\e_{ki}e_{il}\right)+(e_{kk}+e_{ll})e_{kl}-
  (\e_{kk}-\e_{ll})\e_{kl};
  \intertext{for $k>l$:}
  \dd(e_{kl}) &
  =2\sum_{l<i<k}\left(e_{ki}\e_{il}-\e_{ki}e_{il}\right)-(\e_{kk}-\e_{ll})e_{kl}+
  (e_{kk}-e_{ll})\e_{kl},\\
  \dd(\e_{kl}) &
  =2\sum_{l<i<k}\left(-e_{ki}e_{il}+\e_{ki}e_{il}\right)-(e_{kk}+e_{ll})e_{kl}+
  (\e_{kk}-\e_{ll})\e_{kl}.
\end{align}
Here we use the basis $e_{ij}, \e_{ij}\in\mathfrak{q}(n)$, where
$e_{ij}=\diag(E_{ij},E_{ij})\in \mathfrak{q}(n)_0$,
$\e_{ij}=\antidiag(E_{ij},E_{ij})\in \mathfrak{q}(n)_1$, and
$E_{ij}$ are the usual matrix units. Note that
$(e_{ij},\e_{ji})=-(\e_{ij},e_{ji})=1$ and all other inner
products of basis vectors are zero. A basis in $\mathfrak{n}_+$ is
$e_{ij}, \e_{ij}$ with $i<j$, a basis in $\mathfrak{n}_-$ is
$e_{ij}, \e_{ij}$ with $i>j$, and a basis in $\mathfrak{h}$ is
$e_{ii}, \e_{ii}$.
\end{ex}

(Notice that for $\mathfrak{q}(1)$
formulas~\eqref{koskobkanakartane} give zero cobracket, so the odd
bialgebra structure given in Example~\ref{q1} is different from
the above.)

\section{Appendix}

\subsection{Cotangent bundles for dual vector bundles}

Here we give a full proof of Theorem~\ref{testar}. The following
statement, for purely even vector bundles, appeared for the first
time, to my knowledge, in the paper by Mackenzie and
Xu~\cite{mackenzie:bialg}. (Earlier
Tulczyjew~\cite{tulczyjew:1977} considered the particular case of
tangent bundles.)

\begin{thm}\label{cot}
For an arbitrary vector bundle $E\to M$ there is a natural
diffeomorphism of the cotangent bundles $F{\co} T^*E\cong T^*E^*$
which preserves the symplectic structure. It is a natural
transformation in the categorical sense, and $F^2{\co} T^*E\to
T^*E^{**}\cong T^*E$ is the transformation induced by the
multiplication by $-1$ in the fibers of $E$. (The diffeomorphism
$F$ is given by the explicit formula~\eqref{iso0} below.)
\end{thm}
\begin{proof}
Let us denote by $x^a$ local coordinates on $M$, by $y^i$
coordinates in the fiber of $E\to M$, by $y_i$ the corresponding
coordinates in the fiber of the dual bundle $E^*\to M$. Then a
change of coordinates in $E$   has the form
\begin{equation}
  \left\{\begin{aligned}
            x^a&=x^a(x')\\
            y^i&=y^{i'}T_{i'}{}^i(x')
         \end{aligned}\right.
  \quad
\left\{\begin{aligned}
            x^a&=x^a(x')\\
            y_i&=T_{i}{}^{i'}(x')\,y_{i'}
         \end{aligned}\right.
\end{equation}
(Here $T_{i}{}^{i'}T_{i'}{}^j=\delta_i{}^j$.) Denote the conjugate momenta in $T^*E$ and $T^*E^*$
by $p_a, p_i$ and $p_a, p^i$, respectively. For $T^*E$ we obtain
\begin{equation}\label{tstare}
  \left\{\begin{aligned}
            p_a&=\der{x^{a'}}{x^a}p_{a'}+\der{y^{i'}}{x^a}p_{i'}=
            \der{x^{a'}}{x^a}p_{a'}+(-1)^{\at\itt}y^{k'}T_{k'}{}^{i}\der{T_i{}^{i'}}{x^a}p_{i'}\\
            p_i&=\der{y^{i'}}{y^i}p_{i'}=T_{i}{}^{i'}(x')p_{i'}
         \end{aligned}\right.
\end{equation}
and for $T^*E^*$ we obtain
\begin{equation}
  \left\{\begin{aligned}
            p_a&=\der{x^{a'}}{x^a}p_{a'}+\der{y_{k'}}{x^a}p^{k'}=
            \der{x^{a'}}{x^a}p_{a'}+\der{T_{k'}{}^{i}}{x^a}T_{i}{}^{i'}y_{i'}p^{k'}\\
            p^i&=\der{y_{i'}}{y_i}p^{i'}=(-1)^{\itt(\itt+\itt')}T_{i'}{}^{i}(x')p^{i'}.
         \end{aligned}\right.
\end{equation}
(We assume that $\tilde x^a=\tilde a$, $\tilde y^i=\itt$.) Notice
that
$(-1)^{\itt(\itt+\itt')}T_{i'}{}^{i}p^{i'}=(-1)^{\itt+\itt'}p^{i'}T_{i'}{}^{i}$
and
\begin{equation}\label{dtt}
\der{T_{k'}{}^{i}}{x^a}T_{i}{}^{i'}+ (-1)^{\at(\kt'+\itt)}T_{k'}{}^{i}\der{T_{i}{}^{i'}}{x^a}=0,
\end{equation}
hence
$$
\der{T_{k'}{}^{i}}{x^a}T_{i}{}^{i'}y_{i'}p^{k'}=-(-1)^{\at(\kt'+\itt)+\kt'(\kt'+\at)}p^{k'}
T_{k'}{}^{i}\der{T_{i}{}^{i'}}{x^a}y_{i'}= -(-1)^{\at\itt+\kt'}p^{k'}
T_{k'}{}^{i}\der{T_{i}{}^{i'}}{x^a}y_{i'}.
$$
We see that it is possible to rewrite the change of coordinates on $T^*E^*$ as
\begin{equation}\label{cotanestar}
  \left\{\begin{aligned}
            x^a&=x^a(x')\\
            y_i&=T_{i}{}^{i'}(x')\,y_{i'}\\
            p_a&=\der{x^{a'}}{x^a}p_{a'}-(-1)^{\at\itt+\kt'}p^{k'}
            T_{k'}{}^{i}\der{T_{i}{}^{i'}}{x^a}y_{i'}\\
            (-1)^{\itt}p^i&=(-1)^{\itt'}p^{i'}T_{i'}{}^{i}
         \end{aligned}\right.
\end{equation}
which is identical with that for $T^*E$ if we substitute $y_i:=p_i$,
$p^i:=(-1)^{\itt+1}y^i$. Hence, we define the desired  diffeomorphism $F{\co} T^*E\to
T^*E^*$ as
\begin{equation}\label{iso0}
  (x^a,y^i,p_a,p_i)\mapsto (x^a,y_i,p_a,p^i)=(x^a,p_i,p_a,-(-1)^{\itt}y^i).
\end{equation}
Consider the symplectic form $dp_adx^a+dp^idy_i$ on $T^*E^*$.
Substituting~\eqref{iso0}, we get
$dp_adx^a+(-1)^{\itt+1}dy^idp_i=dp_adx^a+dp_idy^i$, which is the
symplectic form on $T^*E$. Thus, the diffeomorphism~\eqref{iso0}
is a canonical transformation. Notice, finally, that
formulas~\eqref{iso0} define the transformation $F$ for all vector
bundles. One only has to be careful with the distinction of left
and right coordinates: our formulas make use of left coordinates
in $E$ and the corresponding right coordinates in $E^*$. (If we
want to change, say, a left coordinate $y^i$ into the right
coordinate, we have to multiply both $y^i$ and the corresponding
momentum by $(-1)^{\itt}$.) Rewritten in left coordinates only,
the formula for $F$ becomes
\begin{equation}\label{iso1}
  F_E{\co}  (x^a,y^i,p_a,p_i)\mapsto (x^a,\bar y_i,p_a,\bar p^i)=(x^a,(-1)^{\itt}p_i,p_a,-y^i)
\end{equation}
where bar is used to denote the left coordinates in $E^*$ and the respective momenta. (The natural
pairing of $E$ and $E^*$ is $y^iy_i=(-1)^{\itt}y^i\bar y_i$.) Hence, for $F_{E^*}$ we will get
\begin{equation}\label{iso2}
  F_{E^*}{\co}  (x^a,\bar y_i,p_a,\bar p^i)\mapsto (x^a,w^i,p_a,q_i)=(x^a,(-1)^{\itt}\bar p^i,p_a,-\bar y_i).
\end{equation}
At the right hand side $w^i$ stand  for  (left) coordinates in
$E^{**}$ and $q_i$ for the corresponding momenta. Notice that the
natural isomorphism between $E^{**}$ and $E$  is, in  left
coordinates, $I_E{\co} (x^a,w^i)\mapsto
(x^a,y^i)=(x^a,(-1)^{\itt}w^i)$. Hence
\begin{equation}\label{iso3}
  I_E\circ F_{E^*}{\co}
  (x^a,\bar y_i,p_a,\bar p^i)\mapsto (x^a,y^i,p_a,p_i)=(x^a,\bar p^i,p_a,-(-1)^{\itt}\bar y_i).
\end{equation}
Taking composition with $F_E$ and   using formula~\eqref{iso1} we
finally obtain:
\begin{equation}\label{iso4}
  I_E\circ F_{E^*}\circ F_{E}{\co}
  (x^a,y^i,p_a,p_i)\mapsto (x^a,y^i,p_a,p_i)=(x^a,-y^i,p_a,-p_i),
\end{equation}
which concludes the proof, if we identify $E^{**}$ with $E$ with the help of $I_E$.
\end{proof}

\subsection{Analog for anticotangent bundles}
\begin{thm} \label{pitestar}
For an arbitrary vector bundle $E\to M$ there is a natural
diffeomorphism (in the categorical sense) of the anticotangent
bundles
\begin{equation}
  F{\co} \Pi T^*E\cong \Pi T^*(\Pi E^*),
\end{equation}
which preserves the canonical odd symplectic structure. It is
given by formula~\eqref{piso0} below. The square $F^2{\co} \Pi
T^*E\to \Pi T^*(\Pi (\Pi E^{*})^{*})\cong \Pi T^*E$ is the
transformation induced by the multiplication by $-1$ in the fibers
of $E$.
\end{thm}
\begin{proof} Let  coordinates in $E$ be denoted, as above, by $x^a$ and $y^i$. Let $x^*_a$
and $y^*_i$ be the corresponding odd momenta. Then, similarly to~\eqref{tstare}, we have
\begin{equation}
  \left\{
  \begin{aligned}
    x^a&=x^a(x')\\
    y^i&=y^{i'}T_{i'}{}^i(x')\\
    x^*_a&=\der{x^{a'}}{x^a}x^*_{a'}+(-1)^{\at\itt}y^{k'}T_{k'}{}^{i}\der{T_i{}^{i'}}{x^a}y^*_{i'}\\
    y^*_i&=T_{i}{}^{i'}(x')x^*_{i'}
  \end{aligned}\right.
\end{equation}
Consider now the antidual bundle $\Pi E^*$. Let $\h_i$ denote the fiber coordinates in $\Pi
E^*$ that are right-contragredient to $y^i$, i.e., the form $y^i\h_i$ gives the invariant
canonical odd pairing of $E$ and $\Pi E^*$. Notice that $\tilde \h_i=\tilde y^i+1=\itt+1$.
Then
\begin{equation}
  \h_i=T_i{}^{i'}(x')\h_{i'}=(-1)^{(\itt'+1)\itt}\h_{i'}T_i{}^{i'}(x'),
\end{equation}
and we get
\begin{equation}
  \left\{
  \begin{aligned}
    x^a&=x^a(x')\\
    \h_i&=T_i{}^{i'}(x')\h_{i'}\\
    x^*_a&=\der{x^{a'}}{x^a}x^*_{a'}+\der{T_{k'}{}^{i}}{x^a}T_{i}{}^{i'}\h_{i'}\h^{*k'} \\
    \h^{*i}&=(-1)^{(\itt+1)\itt'}T_{i'}{}^{i}(x')\h^{*i'}
  \end{aligned}\right.
\end{equation}
for the changes of coordinates  in $\Pi T^*(\Pi E^*)$. Notice that
$\tilde \h^{*i}=\tilde \h_i+1=\itt$. Using formula~\eqref{dtt} and
moving $\h^{*k'}$ to the left, we can rewrite it as
\begin{equation}
  \left\{
  \begin{aligned}
    x^a&=x^a(x')\\
    \h_i&=T_i{}^{i'}(x')\h_{i'}\\
    x^*_a&=\der{x^{a'}}{x^a}x^*_{a'}-(-1)^{\at\itt}\h^{*k'}T_{k'}{}^{i}\der{T_{i}{}^{i'}}{x^a}\h_{i'} \\
    \h^{*i}&=\h^{*i'}T_{i'}{}^{i}(x')
  \end{aligned}\right.
\end{equation}
Hence, we can define the desired diffeomorphism between $\Pi T^*E$ and $\Pi T^*(\Pi E^*)$ by
the formula
\begin{equation}\label{piso0}
  F{\co}  (x^a,y^i,x^*_a,y^*_i)\mapsto (x^a,\h_i,x^*_a,\h^{*i})=(x^a,y^*_i,x^*_a,-y^i).
\end{equation}
Check the   symplectic form: on $\Pi T^*(\Pi E^*)$ it is
$\o=(-1)^{\at+1}dx^*_adx^a+(-1)^{\itt}d\h^{*i}d\h_i$, and
substituting $\h_i:=y^*_i$, $\h^{*i}:=-y^i$ we get
$(-1)^{\at+1}dx^*_adx^a-(-1)^{\itt}dy^idy^*_i=(-1)^{\at+1}dx^*_adx^a+(-1)^{\itt+1}dy^*_idy^i$,
which is exactly the odd symplectic form on $\Pi T^*E$. We omit
the direct check for $F^2$, which is completely similar to that in
the proof of Theorem~\ref{cot}.
\end{proof}

\begin{ex}\label{linbivector}
Let $\mathfrak g$ be a vector space.  Consider an even linear bivector field $P=(1/2)\x^k
P^{ij}_k \x^*_j\x^*_i$ on $\Pi \mathfrak g$. According to our theorem, there is a
diffeomorphism $\Pi T^*\Pi \mathfrak g\to \Pi T^*\mathfrak g^*$ (preserving the odd
symplectic structure), which in coordinates is given by $x_i=\x^*_i$, $x^{*i}=-\x^i$.
Hence, we can transform $P$ into the function
$-(1/2)x^{*k}P^{ij}_k x_jx_i=-({1}/{2})(-1)^{\kt+1}P^{ij}_k
x_jx_ix^{*k}$
on $\Pi T^*\mathfrak g^*$, which is an even antivector field on
$\mathfrak g^*$. Recalling $\t{\co}  \Vect (\mathfrak g^*)\to
\fun(\Pi T^*\mathfrak g^*)$, $X=X_k(x)\lder{}{x_k}\mapsto
(-1)^{\Xt}X_k(x)x^{*k}$, we recover an odd vector field $\hat
Q^*=({1}/{2})(-1)^{\kt+1}P^{ij}_k x_jx_i\lder{}{x_k}$. Because the
symplectic structure is preserved, the field $\hat Q^*$  on
$\mathfrak g^*$ is homological if and only if $P$ is a Poisson
tensor.
\end{ex}

\begin{prop}
The Lie bracket in $\Pi \mathfrak g^*$ defined by the field $\hat Q^*$ is given by the
formulas
\begin{equation}
  [\e^i,\e^j]=(-1)^{\itt+1}\e^k P^{ij}_k.
\end{equation}
Hence, the corresponding Lie--Poisson bracket on $\Pi \mathfrak g$
coincides with the Poisson bracket directly given by the tensor
$P$.
\end{prop}

\bigskip\noindent
{\footnotesize \textbf{Acknowledgements.} The earlier version of
this paper grew out of comments on the talks by D.~Roytenberg and
P.~{\v S}evera at the conference \textit{Poisson 2000} in June
2000 (CIRM, Luminy) written for myself. I thank
Y.~Kosmann-Schwarzbach for the invitation and for the wonderful
atmosphere at the conference. This earlier version was presented
at the workshop on quantization at Warwick in July 2000. I thank
John Rawnsley for the invitation.   I would like to mention that I
first got acquainted with the problem of ``Drinfeld's double'' for
Lie bialgebroids in 1997 at Berkeley from  D.~Roytenberg, who then
was a  Ph.D. student of Alan~Weinstein. I witnessed how his
work~\cite{roytenberg:thesis} came into being. It is a pleasure to
mention our discussions. Also, I would like to thank
H.~Khudaverdian, K.~Mackenzie, P.~{\v Severa} and A.~Weinstein for
inspiring discussions of many topics. Jim Stasheff, Kirill
Mackenzie and Dmitry Roytenberg read versions of this manuscript
at various stages and made many important remarks. I am deeply
grateful to them.

}



\def\cprime{$'$}

\end{document}